\documentclass[11pt,reqno]{amsart}
\usepackage{amsfonts}
\usepackage{graphicx}
\usepackage[usenames]{color}
\usepackage{amssymb}
\usepackage{amsthm}
\usepackage{amssymb}
\usepackage[dvips]{epsfig}

\def\bleu{\textcolor{NavyBlue}}

\def\rouge{\textcolor{BrickRed}}
\def\rose{\textcolor{Lavender}}
\def\jaune{\textcolor{Goldenrod}}
\def\vert{\textcolor{OliveGreen}}

\def\bleu{\textcolor{blue}}

\def\rouge{\textcolor{red}}
\def\rose{\textcolor{magenta}}
\def\jaune{\textcolor{yellow}}
\def\vert{\textcolor{green}}

\def\carre{\bleu{\linethickness{\unitlength}\line(1,0){1}}}
\def\vcarre{\vert{\linethickness{\unitlength}\line(1,0){1}}}
\def\jcarre{\jaune{\linethickness{\unitlength}\line(1,0){1}}}

\numberwithin{equation}{section}
\newtheorem{thm}{Theorem}[section]
\newtheorem{prop}[thm]{Proposition}
\newtheorem{lem}[thm]{Lemma}

\newtheorem{cor}[thm]{Corollary}

\theoremstyle{remark}

\def\charac{{\raise 2pt\hbox{$\chi$}}}
\newcommand{\x}{{\bf x}}
\newcommand{\y}{{\bf y}}
\newcommand{\z}{{\bf z}}

\newcommand{\QQ}{\mathbb{Q}}
\newcommand{\NN}{\mathbb{N}}
\newcommand{\ZZ}{\mathbb{Z}}

\newcommand{\be}{{\beta}}
\newcommand{\s}{{\sigma}}
\newcommand{\de}{{\delta}}
\newcommand{\la}{{\lambda}}

\newcommand{\maj}{\operatorname{maj}}

\newcommand{\Des}{\operatorname{Des}}
\newcommand{\fmaj}{\operatorname{fmaj}}

\newcommand{\nneg}{\operatorname{neg}}

\newcommand{\Sx}{\operatorname{Sch}}
\newcommand{\sx}{\operatorname{s}}
\newcommand{\Hor}{\operatorname{Horiz}}
\newcommand{\Ver}{\operatorname{Vert}}

\newcommand{\Ris}{\operatorname{Ris}}

\title[$B_n$ Action on Coinvariants of $B_n\times B_n$]{Tensorial square of the Hyperoctahedral group Coinvariant Space}
\author[F. Bergeron]{Fran\c cois Bergeron}
\address[F. Bergeron]{D\'epartement de Math\'ematiques\\ Universit\'e
  du Qu\'ebec \`a Mon\-tr\'eal\\ Mont\-r\'eal, Qu\'ebec, H3C 3P8, CANADA}
\email{bergeron.francois@uqam.ca}
\author[R. Biagioli]{Riccardo Biagioli}
\address[R. Biagioli]{LaCIM\\ Universit\'e
  du Qu\'ebec \`a Mont\-r\'eal\\ Montr\'eal, Qu\'ebec, H3C 3P8, CANADA}
\email{biagioli@math.uqam.ca}
\date{\today}
\thanks{F. Bergeron is supported in part by NSERC-Canada and
  FQRNT-Qu\'ebec.}

\begin{document}

\begin{abstract}The purpose of this paper is to give an explicit description of the trivial and alternating components of the irreducible representation decomposition of the bigraded module obtained as the tensor square of the coinvariant space for hyperoctahedral groups.  
\end{abstract}
\maketitle
{ \parskip=0pt\footnotesize \tableofcontents}
\parskip=8pt  

\section{Introduction.} \label{intro}
The purpose of this paper is to study the so called ``diagonal action'' of a the hyperoctahedral group on the tensor square of its coinvariant space. One of the many reasons to study this space is that it contains the space of ``diagonal coinvariants'' of $B_n$ in a very natural way (see e.g., \cite{bergeron_lamonta}, \cite{gordon},\cite{haimanhilb}, and \cite{lamontagne}). 

The paper is organized as follows. We start off with a rapid survey of classical results regarding coinvariant spaces of finite reflection groups, followed with the implications about the tensor square of these same spaces. We then specialize our discussion to hyperoctahedral groups, recalling in the process the main aspects of their representation theory.
 
\section{Reflection group action on the polynomial ring}\label{weyl_action}

For any finite reflection group $W$, on a finite dimensional vector space $V$ over $\QQ$, there corresponds a natural action of $W$ on the polynomial ring $\QQ[V]$. In particular, if $\x=x_{1},\ldots ,x_{n}$ is a basis of $V$ then $\QQ[V]$ can be identified with the ring  $\QQ[\x]$ of polynomials in the variables $x_{1},\ldots ,x_{n}$. As usual, we denote
\begin{equation*}
     w\cdot p(\x)=p(w\cdot \x),
\end{equation*}
the action in question. It is clear that this action of $W$ is degree preserving, thus making natural the following considerations.  Let us denote $\pi_d(p(\x))$ the degree $d$ homogeneous  component of a polynomial $p(\x)$.  The ring  $Q:=\QQ[\x]$ is graded by degree, hence
 \begin{equation*}
     Q\simeq\bigoplus_{d\geq 0} Q_d,
  \end{equation*}
where $Q_d:=\pi_d(Q)$ is the {\em degree $d$ homogeneous component} of $Q$. Recall that a subspace $S$ is said to be {\em homogeneous} if $\pi_d(S)\subseteq S$ for all $d$. Whenever this is the case, we clearly have $ S=\bigoplus_{d\geq 0} S_d$, with $S_d:=S\cap Q_d$, and thus it makes sense to consider the {\em Hilbert series} of $S$:
\begin{equation*}
     H_q(S):=\sum_{m \geq 0}\dim(S_d) \; q^d.
\end{equation*} 
The motivation, behind the introduction of this formal power series in $q$, is that it condenses in a efficient and compact form the information for the dimensions of each of the $S_d$'s. To illustrate, it is not hard to show that the Hilbert series of $Q$ is
 simply
 \begin{equation}\label{hilbert_Q}
     H_q(Q)=\frac{1}{(1-q)^n},
\end{equation}
which is equivalent to the infinite list of statements
 \begin{equation*}
     \dim(Q_d)={n+d-1 \choose d},
\end{equation*}
one for each value of $d$.

We will be particularly interested in {\em invariant subspaces} $S$ of $Q$, namely those for which $w\cdot S\subseteq S$, for all $w$ in $W$.  Clearly whenever $S$ is homogeneous, on top of being invariant, then each of these homogeneous component, $S_d$, is also a $W$-invariant subspace. One important example of homogeneous invariant subspace, denoted $Q^{W}$, is the set of {\em invariant polynomials}. These are the polynomials $p(\x)$ such that
\begin{equation*}
     w\cdot p(\x)= p(\x).
 \end{equation*}
 Here, not only the subspace $Q^W$ is invariant, but all of its elements are.
It is well know that $Q^{W}$ is in fact a subring of $Q$, for which one can find generator sets of $n$ homogeneous algebraically independent elements, say $f_1,\ldots, f_n$, whose respective degrees will be denoted $d_1,\ldots,d_n$. Although the $f_i$'s are not uniquely characterized, the $d_i$'s are basic numerical invariants of the group, called the {\em degrees of $W$}. Any $n$-set $\{f_1,\ldots, f_n\}$ of invariants with these properties is called a {\em set of basic invariants} for $W$. It follows that the Hilbert series of $Q^{W}$ takes the form:
\begin{equation}\label{hilbert_invariant}
   H_q(Q^W)=\prod_{i=1}^n\frac{1}{1-q^{d_i}}.
\end{equation}
Now, let $\mathcal{I}_W$ be the ideal of $Q$ generated by constant term free elements of $Q^{W}$. The {\em coinvariant space} of $W$ is defined to be
\begin{equation}\label{coinvariant}
     Q_W:=Q/\mathcal{I}_W.
 \end{equation}
Observe that, since $\mathcal{I}_W$ is an homogeneous subspace of $Q$, it follows that the ring $Q_W$ is naturally graded by degree. Moreover,  $\mathcal{I}_W$ being $W$-invariant, the group $W$ acts naturally on $Q_W$. In fact, it can be shown that $Q_W$ is actually isomorphic to the \emph{left regular representation} of $W$ (For more on this see \cite{orbit}
or \cite{stanley}). It follows that the dimension of $Q_W$ is exactly the order of the group $W$. We can get a finer description of this fact using a theorem of Chevalley (see \cite[Section 3.5]{humphreys}) that can be stated as follows. There exists a natural isomorphism of $\QQ W$-module\footnote{Here, as usual, $\QQ W$ stands for the group algebra of $W$, and the term {\em module} underlines that we are extending the action of $W$ to its group algebra.}:
\begin{equation}\label{chevalley}
   Q\simeq Q^W \otimes \; Q_W.
\end{equation}
We will use strongly this decomposition in the rest of the paper. One immediate consequence, in view of (\ref{hilbert_Q}) and (\ref{hilbert_invariant}), is that
\begin{eqnarray}\label{hilbert_Q_W}
  H_q(Q_W)&=& \prod_{i=1}^n \frac{1-q^{d_i}}{1-q}\\
       &=& (1+\cdots + q^{d_1-1}) \cdots (1+\cdots + q^{d_n-1}).\nonumber
\end{eqnarray}
We now introduce another important $\QQ W$-module for our discussion. To describe it, let us first introduce a $W$-invariant scalar product on $Q$, namely
\begin{equation*}
  \langle p,q \rangle:=p(\partial\x)q(\x)|_{\x=0}.
\end{equation*}
Here $p(\partial\x)$ stands for the linear operator obtained by replacing each variable $x_i$, in the polynomial $p(\x)$, by the partial derivative $\partial x_i$ with respect to $x_i$. We have denoted above by $\x=0$ the simultaneous substitutions $x_i=0$, one for each $i$. With this in mind, we define  the {\em space of $W$-harmonic polynomials}:
\begin{equation}\label{harmonic}
     \mathcal{H}_W:=\mathcal{I}_W^{\perp}
\end{equation}
where, as usual, $\perp$ stands for orthogonal complement with respect to the underlying scalar product. Equivalently, since $\mathcal{I}_W$ is can be described as the ideal generated (as above) by a basic set $\{f_1,\ldots , f_n\}$ of invariants, then a polynomial $p(\x)$ is in $ \mathcal{H}_W$ if and only if
\begin{equation*}
      f_k(\partial x) p(\x)=0, \qquad \hbox{for all}\ k\geq 0.
 \end{equation*}
It can be shown that the spaces $\mathcal{H}_W$ and $Q_W$ are actually isomorphic as graded $\QQ W$-modules \cite{steinberg}. On the other hand, it is easy to observe that $\mathcal{H}_W$ is closed under partial derivatives.
These observations, together with a further remark about characterizations of reflection groups contained in Chevalley's Theorem,  make possible an explicit description of $\mathcal{H}_W$ in term of the {\em Jacobian determinant}:\begin{equation}\label{vandermonde}
      \Delta_W(\x):=\frac{1}{|W|}\det \frac{\partial f_i}{\partial x_j},
\end{equation} 
where the $f_i$'s form a set of basic $W$-invariants.  This polynomial is also simply denoted $\Delta(\x)$, when the underlying group is clear. One can show that this polynomial is well defined (up to a scalar multiple) in that it does not depend on the actual choice of the $f_i$'s (see \cite[Section 3.13]{humphreys}). It can also be shown that $\Delta_W$ is the unique (up to scalar multiple) $W$-harmonic polynomial of maximal degree, and we have
\begin{equation}\label{harmonic_explicit}
    \mathcal{H}_W=\mathcal{L}_{\partial}[\Delta_W(\x)],
\end{equation}
where $\mathcal{L}_{\partial}$ stands as short hand for ``linear span of all partial derivatives of''.  Another important property of $\Delta=\Delta_W$ is that it allows an explicit characterization of all {\em $W$-alternating polynomials}. Recall that, $p(\x)$ is said to be $W$-alternating if and only if 
    $$w\cdot p(\x) = \det(w)\, p(\x),$$
where, to make sense out of $\det(w)$, one interprets $w$ as linear transformation.
The pertinent statement is that $p(\x)$ is alternating if and only if it can be written as
   $$p(\x)= f(\x)\, \Delta(\x),$$
with $f(\x)$ in $Q^W$.  In other words, $\Delta$ is the minimal $W$-alternating polynomial.
Thus, in view of (\ref{hilbert_invariant}) and (\ref{vandermonde}), the Hilbert series of the homogeneous invariant subspace $Q^{\pm}$, of $W$-alternating polynomials,  is simply
\begin{equation}\label{hilbert_alternating}
   H_q(Q^{\pm})=\prod_{i=1}^n\frac{q^{d_i-1}}{1-q^{d_i}}.
\end{equation}
The point of all this is that we can reformulate the decomposition given in (\ref{chevalley}) as
 \begin{equation}\label{chev_harm}
     Q\simeq Q^W \otimes \mathcal{H}_W,
  \end{equation}
  with both $Q^W$ and $\mathcal{H}_W$ $W$-submodules of $Q$.  In other words, there is a unique decomposition of any polynomial $p(\x)$ of the form
     $$p(\x)=\sum_{w\in W} f_w(\x)\, b_w(\x),$$
 for any given basis $\{ b_w(\x)\ |\ w\in W\}$ of $\mathcal{H}_W$, with the $f_w(\x)$'s invariant polynomials. Recall here that $\mathcal{H}_W$ has dimension equal to $|W|$. As we will see in particular instances, there are natural choices for such a basis.

\section{Diagonally invariant and alternating polynomials}\label{diagonallypolys}
We now extend our discussion to the ring 
    $$R=\QQ[\x,\y]:=\QQ[x_1,\ldots,x_n,y_1,\ldots,y_n],$$ 
of polynomials in two sets of $n$ variables, on which we want to study the {\em diagonal action} of $W$, namely such that:
\begin{equation}\label{diagonal}
    w\cdot p(\x,\y)=p(w\cdot\x,w\cdot\y),
\end{equation}
for $w\in W$. In this case, $W$ does not act as a reflection group on the vector space $\mathcal{V}$ spanned by the $x_i$'s and $y_j$'s, so that we are truly in front of a new situation, as we will see in more details below. By comparison, the results of Section \ref{weyl_action} would still apply to $R$ if we would rather consider the action of $W\times W$,  for which
\begin{equation}\label{tensorial}
   (w,\tau)\cdot p(\x,\y)=p(w\cdot\x,\tau\cdot\y),
\end{equation}
when $(w,\tau)\in W \times W$, and $p(\x,\y)\in R$. Indeed, this does correspond to an action of $W\times W$ as a reflection group on $\mathcal{V}$. Each of these two contexts give rise to a notion of invariant polynomials in the same space $R$. Notation wise, we naturally distinguish these two notions as follows. On one hand we have the subring $R^W$ of {\em diagonally invariant polynomials}, namely those for which
    $$p(w\cdot \x, w\cdot \y)=p(\x,\y);$$
and, on the other hand, we get the subring $R^{W\times W}$, of invariants polynomials of the tensor action (\ref{tensorial}), as a special case of the results described in Section \ref{weyl_action}. Observe that
\begin{equation}\label{tensordecomposition}
    R^{W\times W} \simeq \QQ[\x]^{W} \otimes \QQ[\y]^{W}.
\end{equation}
In view of this observation, we will called $R^{W\times W}$ the  {\em tensor invariant algebra}. It is easy to see that  $R^{W\times W}$ is a subring of $R^W$.

The ring $R$ is naturally ``bigraded'' with respect to ``bidegree''. To make sense out of this, let us recall the usual vectorial notation for monomials:
     $$\x^{\bf a}\y^{\bf b}:=x_1^{a_1}\cdots x_n^{a_n}\,y_1^{b_1}\cdots y_n^{b_n},$$
with ${\bf a}=(a_1,\ldots,a_n)$ and ${\bf b}=(b_1,\ldots, b_n)$ both in $\NN^n$. Then, the {\em bidegree} of $\x^{\bf a}\y^{\bf b}$ is simply $( |{\bf a}|, |{\bf b}|)$, where $|{\bf a}|$ stands for the sum of the components of ${\bf a}$, and likewise for ${\bf b}$. If we now introduce the linear operator $\pi_{k,j}$ such that
    $$\pi_{k,j}(\x^{\bf a}\y^{\bf b}):=\begin{cases}
      \x^{\bf a}\y^{\bf b}& \text{if}\  |{\bf a}|=k\ {\rm and}\  |{\bf b}|=j, \\
      0 & \text{otherwise},
\end{cases}$$
then a polynomial $p(\x,\y)$ is said to be {\em bihomogeneous} of bidegree $(k,j)$ if and only if
  $$\pi_{k,j}(p(\x,\y))=p(\x,\y).$$
The notion of {\em bigrading} is then obvious. For instance, we have the bigraded decomposition: 
\begin{equation*}
   R=\bigoplus_{k,j} R_{k,j},
\end{equation*}
with $R_{k,j}:=\pi_{k,j}(R)$. Naturally, a subspace $S$ of $R$ is said to be {\em bihomogeneous} if $\pi_{k,j}(S)\subseteq S$ for all $k$ and $j$, and it is just as natural to consider the {\em bigraded Hilbert series}:
   $$H_{q,t}(S):=\sum_{k,j} \dim(S_{k,j}) q^k t^j,$$
with $S_{k,j}:=S\cap R_{k,j}$. Since it is clear that
\begin{equation*}
    R\simeq \QQ[\x]\otimes \QQ[\y],
 \end{equation*}
from (\ref{hilbert_Q}) we easily get
\begin{equation}\label{hilbert_Rxy}
    H_{q,t}(R)=\frac{1}{(1-q)^n}\,\frac{1}{(1-t)^n}.
\end{equation}
Furthermore, in view of (\ref{tensordecomposition}) and (\ref{hilbert_invariant}), the bigraded Hilbert series of $R^{W\times W}$ is simply
\begin{equation}\label{hilb-RWW}
     H_{q,t}(R^{W\times W})= \prod_{i=1}^{n}\frac{1}{(1-q^{d_i})(1-t^{d_i})}.
\end{equation}
Let $R_{W\times W}$ and $\mathcal{H}_{W\times W}$ the spaces of coinvariants and harmonics of $W \times W$, defined in (\ref{coinvariant}) and (\ref{harmonic}), respectively. From (\ref{hilbert_Rxy}) and (\ref{hilb-RWW}), we conclude that
\begin{eqnarray}\label{harmonic_tens}
   H_{q,t}(R_{W\times W})&=&H_{q,t}(\mathcal{H}_{W\times W})\nonumber \\
      &=&  \prod_{i=1}^n \frac{1-q^{d_i}}{1-q} \prod_{i=1}^n \frac{1-t^{d_i}}{1-t}.
 \end{eqnarray}
In fact, we have $(W\times W)$-module isomorphisms of bigraded spaces
\begin{eqnarray}
    \mathcal{C} &\simeq& Q_W\otimes Q_W\\
    \mathcal{H} &\simeq&  \mathcal{H}_W\otimes  \mathcal{H}_W. \label{d-H}
\end{eqnarray}
Here, and from now on, we simply denote $\mathcal{C}$ the space of coinvariants of $W\times W$, and $\mathcal{H}$ the space of harmonics of $W\times W$. Recalling our previous general discussion, the spaces $\mathcal{C}$ and $\mathcal{H}$ are isomorphic as bigraded $W$-modules. Summing up, and considering $W$ as a {\em diagonal} subgroup of $W\times W$ (i.e.: $w\mapsto (w,w)$) we get an isomorphism of $W$-module
 \begin{equation}\label{the_iso}
     R \simeq \QQ[\x]^W\otimes \QQ[\y]^W\otimes \mathcal{H}
 \end{equation}
 from which we deduce, in particular, that
 \begin{equation}\label{the_iso_inv}
     R^W \simeq \QQ[\x]^W\otimes \QQ[\y]^W\otimes \mathcal{H}^W,
 \end{equation}
 where
    $$ \mathcal{H}^W:= R^W\cap  \mathcal{H}.$$
Similarly, for the $W$-module of {\em diagonally alternating polynomials}
\begin{equation}
    R^{\pm}:=\{p(\x,\y) \in R \ |\  p(w\cdot \x,w\cdot \y)=\det(w)\, p(\x,\y)\},
\end{equation}
we have the decomposition:
 \begin{equation}\label{the_iso_alt}
    R^{\pm} \simeq \QQ[\x]^W\otimes \QQ[\y]^W\otimes \mathcal{H}^{\pm},
 \end{equation}
where
    $$   \mathcal{H}^{\pm}:= R^{\pm}\cap  \mathcal{H}.$$
Thus the two spaces  $\mathcal{H}^W$ and $\mathcal{H}^{\pm}$, respectively of diagonally symmetric and diagonally alternating harmonic polynomials, play a special role in the understanding of of  $R^{W}$ and $R^{\pm}$. As we will see below, they are also very interesting on their own. Clearly, $\mathcal{H}^W \simeq \mathcal{C}^W $ and  $\mathcal{H}^{\pm} \simeq \mathcal{C}^\pm$. Nice combinatorial descriptions of these two last spaces will be given in the case of Weyl groups of type $B$.  

\section{The Hyperoctahedral group $B_n$}\label{hyper_group}
The hyperoctahedral group $B_{n}$ is the group of  {\em
signed permutations}  of the set
$[n]:=\{1,2,\ldots,n\}$. More precisely, it is obtained as the wreath product, $\ZZ_2 \wr S_n$, of the ``sign change'' group $\ZZ_2$ and the symmetric group $S_n$. In one line notation, elements of $B_n$ can be written as 
      $$\be=\be(1)\be(2)\cdots \be(n),$$
 with each $\beta(i)$ an integer whose absolute value lies in $[n]$. Moreover, if we replace in $\beta$ each these $\beta(i)$'s by their absolute value, we get a permutation. We often denote the negative entries with an overline, thus $\overline{2}\,\overline{1}\,\overline{5}\,4\,3 \in B_5$. 

The action of $\beta$ in $B_n$ on polynomials is entirely characterized by its effect on variables:
     $$\beta\cdot x_i= \pm x_{\sigma(i)},$$
 with the sign equal to the sign of $\beta(i)$, and $\sigma(i)$ equal to its absolute value.
The $B_n$-invariant polynomials  are thus simply the usual symmetric polynomials 
in the square of the variables. We will write 
     $$f({\bf x}^2):=f(x_1^2,x_2^2,\ldots, x_n^2).$$
Hence a set of basic invariant for $Q^{B_n}$ is given by the {\em power sum  symmetric polynomials}:
\begin{equation*}
     p_{j}({\bf x}^2)=\sum_{1\leq k \leq n} x_{k}^{2j},
 \end{equation*}
with $j$ going from $1$ to $n$. Thus $2,4,\ldots,2n$ are the degrees of $B_n$ and from (\ref{hilbert_invariant}) we get 
\begin{equation}\label{Hilbertinvariants}
H_q(Q^{B_n})= \prod_{i=1}^{n}\frac{1}{(1-q^{2i})}.
\end{equation}
It also follows that the Jacobian determinant (\ref{vandermonde}) is
\begin{equation*}
    \Delta(\x)=\frac{1}{2^n n!} \det
     \begin{pmatrix}
             1 & 1 &\ldots & 1\\
             2 x_1 & 2x_2 & \ldots &2 x_n\\
              \vdots & \vdots & \ddots & \vdots\\
              2n x_1^{2n-1} & 2n x_2^{2n-1} & \ldots & 2n x_n^{2n-1}
     \end{pmatrix}
\end{equation*}
which  factors out simply as
\begin{equation}\label{jacobian}
    \Delta(\x)=x_1\cdots x_n \prod_{1\leq i<j\leq n} (x_i^2-x_j^2).
\end{equation}
Now, an easy application of Buchberger's criteria (see e.g., \cite{CLO}) shows that the set 
\begin{equation}\label{grobner}
   \{h_k(x_k^2,\ldots,x_n^2) \ |\  1\leq k \leq n \}
\end{equation}
is a Gr\"obner basis for the ideal $\mathcal{I}=\mathcal{I}_{B_n}$. Here, we are using the lexicographic monomial order (with the variables ordered as $x_1>x_2>\ldots>x_n$), and $h_k$ denotes the $k^{\rm th}$ {\em complete homogeneous symmetric polynomial}. It follows from the corresponding theory that a linear basis for the coinvariant space $Q_{B_n}=Q/\mathcal{I}$ is given by the set \begin{equation}\label{basis_coinvariant}
   \{\x^{\bf \epsilon}+\mathcal{I}\ |\  {\bf \epsilon}=(\epsilon_1,\ldots,\epsilon_n),\  {\rm with}\  0 \leq \epsilon_i < 2i\}.
\end{equation}
These monomials are exactly those that are not divisible by any of the leading terms 
   $$x_1^2,\ldots ,x_k^{2k},\ldots, x_n^{2n}$$
of the polynomials in the Gr\"obner basis (\ref{grobner}).
The linear basis (\ref{basis_coinvariant})  is sometimes called the {\em Artin basis} of the coinvariant space. 
If we systematically order the terms of polynomials in decreasing lexicographic order, it is then easy to deduce, from (\ref{harmonic_explicit}) and (\ref{basis_coinvariant}), that the set
\begin{equation*}
    \{\partial \x^{\bf \epsilon}\Delta(\x)\ |\ {\bf \epsilon}=(\epsilon_1,\ldots,\epsilon_n), 0 \leq \epsilon_i < 2i\}.
\end{equation*}
is a basis of the module of $B_n$-harmonic polynomials. This makes it explicit that $2^nn!$ is the dimension of both $Q_{B_n}$ and $\mathcal{H}_{B_n}$. We will often go back and forth between $Q_{B_n}$ and $\mathcal{H}_{B_n}$, using the fact that they are isomorphic as graded representations of $B_n$.

Another basis of the space of coinvariants, called the {\em descent basis}, will be useful for our purpose. Let us first introduce some ``statistics'' on $B_n$ that also have an important role in our presentation. We start by fixing the following linear order on $ \ZZ $: 
\begin{equation*}
\bar{1}\prec \bar{2} \prec \cdots \prec \bar{n}\prec \cdots \prec 0\prec 1\prec 2\prec \cdots \prec n\prec \cdots .
\end{equation*}
Then, following \cite{adin}, we define the \emph{flag-major
index} of $ \beta \in B_{n} $ by 
\begin{equation}\label{fmaj}
    \fmaj(\beta ):=2\maj(\beta )+\nneg(\beta ),
\end{equation}
where $\nneg(\be)$ is just the number of the negative entries in $\beta$, and $\maj(\beta)$ is the usual {\em major index} of an integer sequence, i.e., 
     \[\maj(\be)=\sum_{i \in \Des(\be)} i.\] 
Here, $\Des(\be)$ stands for the {\em descent set} of $\beta$, namely
      $$\Des(\be):=\{i \in [n-1]\ | \  \be_i \succ \be_{i+1}\}.$$
For example, with $\be=\bar{2}\,\bar{1}\,\bar{5}\,4\,3$, we get $\Des(\be)=\{1,4\}$, $\maj(\be)=5$, $\nneg(\be)=3$, and $\fmaj(\be)=15$.
It will be handy to localize these three statistics setting, for
$i \in [n]$:
\begin{eqnarray}
f_{i}(\be)&:=&2d _{i}(\be )+\varepsilon_{i}(\be ),\ {\rm with}\label{fmaj_i}\\
  \varepsilon_{i}(\be)&:=&\begin{cases}
     1 & \textrm{if }\be(i)<0,\ {\rm and}\\
      0 & \text{otherwise},
\end{cases}
  \label{varepsilon}\\
d_i(\be)&:=&\# \{j \in \Des(\be)\ |\ j \geq i \}.\label{desc_i}
\end{eqnarray}
As is shown in \cite{brenti}, the set $\{\x_{\be }+\mathcal{I} \ |\  \be \in B_{n}\},$ with
\[\x_{\be}:=\prod _{i=1}^{n}x_{\sigma (i)}^{f_{i}(\be)},\] 
is another linear basis of the coinvariant space $Q_{B_n}$, if $\sigma(i)$ denotes the absolute value of $\beta(i)$. Note that each monomial $\x_{\be}$ has precisely degree $\fmaj(\be)$ so that, in view of (\ref{hilbert_Q_W}) we get
\begin{equation*}
   \sum_{\beta\in B_n} q^{\fmaj(\beta)}=\prod_{j=1}^n \frac{1-q^{2j}}{1-q}.
\end{equation*}

\section{Plethystic substitution}
To go on with our discussion, it will be particularly efficient to use the notion of ``plethystic
substitution''. Let ${\bf z}=z_1,z_2,z_3,\ldots$  and $\bar{\bf z}=\bar{z}_1,\bar{z}_2,\bar{z}_3,\ldots$ be two infinite sets of ``formal'' variables, and denote $\Lambda(\z)$  (resp. $\Lambda(\bar{\bf z})$) the ring of symmetric functions\footnote{Here, the term ``function'' is used to emphasize that we are dealing with infinitely many variables.}  in these variables ${\bf z}$ (resp. $\bar{\bf z}$). 
It is well known that any classical linear basis of $\Lambda(\z)$ is naturally indexed by partitions. We usually denote $\ell=\ell(\lambda)$ the number of parts of $\lambda\vdash n$ ($\lambda$ ``a partition  of'' $n$). In accordance with the notation of \cite{macdonald},  we further denote
  $$p_\lambda(\z):=p_{\lambda_1}(\z)p_{\lambda_2}(\z)\cdots p_{\lambda_\ell}(\z),$$
the power sum symmetric {\em function} indexed by the partition $\lambda=\lambda_1\lambda_2\cdots \lambda_\ell$. The homogeneous degree $n$ {\em complete} and {\em elementary} symmetric functions are respectively denoted by $h_n(\z)$ and $e_n(\z)$. Recall that we have
 \begin{eqnarray*}
     \sum_{n\geq 0} h_n(\z) t^n &=&\exp\Big(\sum_{k\geq 1} p_k(\z)t^k/k\Big),\quad {\rm and}\\
     \sum_{n\geq 0} e_n(\z) t^n &=&\exp\Big(\sum_{k\geq 1}(-1)^{k-1} p_k(\z)t^k/k\Big).
 \end{eqnarray*}
Our intent here is to use symmetric functions expressions, obtained by ``plethystic substitution'', to encode characters.  A {\em plethystic substitution} $u[\mathbf{w}]$, of an expression $\mathbf{w}$ into a symmetric function $u$, is defined as follows. The first ingredient used to ``compute'' the resulting expression, is the fact that such a substitution is both additive and multiplicative:
\begin{eqnarray*}
    (u+v)[\mathbf{w}]&=&u[\mathbf{w}]+v[\mathbf{w}] \\
    (uv)[\mathbf{w}]&=&u[\mathbf{w}]\,v[\mathbf{w}].
\end{eqnarray*}    
Moreover, a plethystic substitution into a power sum $p_k$ is defined to result in replacing all variables in $\mathbf{w}$ by their $k^{\rm th}$ power. In particular, this makes such a substitution linear in the argument
   $p_k[\mathbf{w}_1+\mathbf{w}_2]=p_k[\mathbf{w}_1]+p_k[\mathbf{w}_2].$
Summing up, we get
   $$u[\mathbf{w}]:=\sum_{\mu} a_\mu \prod_{i=1}^{\ell(\mu)} p_{\mu_i}[\mathbf{w}]$$
whenever the expansion of $u$ in power sum is $u(\z):=\sum_{\mu} a_\mu p_{\mu}(\z)$.

A further useful convention, in this context, is to denote sets of variables $z_1,z_2,z_3,\ldots$, as (formal) sums $\z=z_1+z_2+z_3+\ldots$.
 This has the nice feature that the result of the plethystic substitution of a ``set'' of variables into a power sum $p_k$:
  \begin{eqnarray*}
     p_k[\z] &=& p_k[z_1+z_2+z_3+\ldots]\\
                &=&  z_1^k+z_2^k+z_3^k+\ldots  
   \end{eqnarray*}
is exactly the corresponding power sum in this set of variables. As an illustration of the many useful formulas one can get using this approach, we have
  \begin{equation}\label{hn_sum}
   h_n[\mathbf{w}_1+\mathbf{w}_2] = \sum_{k=0}^n h_k[\mathbf{w}_1] h_{n-k}[\mathbf{w}_2]
\end{equation}
as well as
\begin{equation}\label{cauchy}
   h_n[\mathbf{w}_1\mathbf{w}_2] = \sum_{\lambda\vdash n} s_\lambda[\mathbf{w}_1] s_\lambda[\mathbf{w}_2],
\end{equation}
where the $s_\lambda$'s are the classical {\em Schur functions}. Recall that
\begin{equation}\label{table_charac}
   s_\lambda(\z)=\sum_{\mu} \charac_\mu^\lambda\,\frac{p_\mu(\z)}{z_\mu},
\end{equation}      
where $\charac_\mu^\lambda$ is the value at $\mu$ of the irreducible character of $S_n$ associated to a partition $\lambda$, and we have set
\begin{equation}\label{def-z}
z_{\mu}:=1^{k_1}k_1!\,2^{k_2}k_2!\cdots n^{k_n}k_n!
\end{equation}
whenever $\mu$ has $k_i$ parts of size $i$. It is well known that conjugacy classes of $S_n$ are naturally indexed by partitions, hence (\ref{table_charac}) is an encoding of the character table of $S_n$. Moreover, the $s_\lambda$'s have a natural role in computations regarding characters of $S_n$, through the Frobenius characteristic transformation. Recall that this is the symmetric function associated to a representation $\mathcal{V}$ of $S_n$, in the following manner 
\begin{equation}\label{frob_Sn}
   \mathcal{F}(\mathcal{V}):=\frac{1}{n!}\sum_{\sigma\in S_n} \charac^{\mathcal{V}} (\sigma) p_{\mu(\sigma)},
\end{equation}
where $\mu(\sigma)$ is the partition giving the cycle index structure of $\sigma$.
It is interesting to illustrate this with the regular representation of $S_n$, for which the Frobenius (characteristic) is
\begin{equation}\label{regular_Sn}
    p_1(\z)^n=\sum_{\lambda \vdash n} f_\lambda s_\lambda(\z).
 \end{equation}
The left hand side corresponds to the direct computation of (\ref{frob_Sn}), and the right hand side describe the usual decomposition of the regular representation in terms of irreducible representations. Here is exemplified the fact that the $s_\lambda$'s correspond to irreducible representations of $S_n$, through the Frobenius characteristic.  A well known fact of representation theory says that the coefficients $f_\lambda$, in (\ref{regular_Sn}), are both the dimension of irreducible representations of $S_n$, and the multiplicities of these in the regular representation. 
As is also very well known, these values are given by the hook length formula.

Formulas  (\ref{hn_sum}) and (\ref{cauchy}) are special cases of the more general formulas:
\begin{eqnarray*}
  s_\lambda[\mathbf{w}_1+\mathbf{w}_2]&=&\sum_{\nu \subseteq \lambda}
         s_\nu[{\bf w}_1] s_{\lambda/\nu}[{\bf w}_2]\\
  s_\lambda[\mathbf{w}_1-\mathbf{w}_2]&=&\sum_{\nu \subseteq \lambda}
         (-1)^{|\lambda/\nu|} s_\nu[{\bf w}_1] s_{\lambda'/\nu'}[{\bf w}_2]\\
 s_{\lambda/\mu}[-\mathbf{w}]&=&
         (-1)^{|\lambda/\mu|} s_{\lambda'/\mu'}[{\bf w}].\\
\end{eqnarray*}
As a particular case of this last formula, we get
  $$h_n[-{\bf w}]= (-1)^n e_n[{\bf w}],$$
since $s_{(n)}=h_n$ and $s_{1^n}=e_n$. 
Moreover, we have 
\begin{equation}\label{h-to-p}
h_n[{\bf w}]= \sum_{\mu \vdash n}\frac{1}{z_{\mu}}p_{\mu}[{\bf w}].
\end{equation}
It will be useful also to introduce 
	$$ \Omega[{\bf w}]:=\sum_{n\geq 0} h_n[{\bf w}].$$
Note that 
\begin{equation}\label{omega}
\Omega[{\bf w}_1+{\bf w}_2]=\Omega[{\bf w}_1]\Omega[{\bf w}_2].
\end{equation} 
 
\section{Frobenius characteristic of $B_n$-modules}
For the purpose of this work, we are actually interested in explicit decompositions of our graded (or bigraded) $B_n$-modules into irreducible representations. As a step toward this end, we wish to compute characters of each (bi)homogeneous of the modules considered.  As we will see, this is best encoded using the ``(bi)graded Frobenius characteristic'' as defined below. But first, let us recall some basic facts about the representation theory of hyperoctahedral groups. Conjugacy classes of the group $ B_{n} $, and thus irreducible characters, are naturally parametrized by ordered pairs $ (\mu^+ ,\mu^- ) $ of partitions such that the total sum of their parts is equal to $ n $ (see, e.g., \cite{MD}). In fact, elements $\beta$ of any given conjugacy class of $B_n$ are exactly characterized by their {\em signed cycle type} 
     $$ \mu(\beta):=(\mu^+ (\beta),\mu^-(\beta) ) $$
where the parts of $ \mu^+(\beta)  $
correspond to sizes of ``positive'' cycles in $\beta$, and the parts of $ \mu^-(\beta) $ correspond to ``negative'' cycles sizes. The {\em sign} of a cycle is simply $(-1)^k$, where $k$ is the number of signed elements in the cycle.

To make our notation in the sequel more compact, we introduce the ``bivariate'' power sum
    $$p_{\mu,\nu}(\z,\bar{\z}):=p_\mu(\z) p_\nu(\bar{\z}).$$
With this short hand notation in mind, the {\em bigraded Frobenius characteristic of type $B$} of an invariant bihomogeneous submodule $S$ of $R$ is then defined to be the series
\begin{equation}\label{bigraded_frobenius}
 \mathcal{F}^B_{q,t}(S):=\sum_{k,j}q^k t^j \mathcal{F}^B(S_{k,j}).
\end{equation}
Here $\mathcal{F}^B$ denotes the {\em Frobenius characteristic of type $B$} of a $B_n$-invariant submodule $T$:
\begin{equation}\label{simple_frobenius}
 \mathcal{F}^B(T):=\sum_{|\mu|+|\nu|=n} 
                     \charac^{T}(\mu,\nu)\, 
                     \frac{p_{\mu,\nu}(\z,\bar{\z})}{z_\mu z_\nu}.
\end{equation}
where $ \charac^{T}(\mu,\nu)$ is the value of the  character of $T$ on elements of the conjugacy class characterized by $(\mu,\nu)$, and $z_{\mu}$ is defined in (\ref{def-z}). Now, it is a well established fact  (see e.g., \cite{Stemb}) that irreducible representations of $B_n$ are  mapped bijectively, by this Frobenius characteristic $\mathcal{F}^B$, to products of Schur functions of the form $s_{\la}[\z+\bar{\z}] s_{\rho}[\z-\bar{\z}]$, with $|\lambda|+|\rho|=n$. In other words, there is a natural indexing by such pairs $(\lambda,\rho)$ of a complete set $\{\mathcal{V}^{\lambda,\rho}\}_{|\lambda|+|\rho|=n}$ of irreducible representations of $B_n$, such that
\begin{equation}\label{stembridge}
    \mathcal{F}^B(\mathcal{V} ^{\la,\rho})=s_{\la}[\z+\bar{\z}] s_{\rho}[\z-\bar{\z}].
\end{equation}
Thus, when expressed in terms of Schur functions, the Frobenius characteristic describes the decomposition into irreducible characters of each bihomogeneous component of a bigraded $B_n$-module $S$. In other words, we  have
\begin{equation}\label{frob_in_schur}
    \mathcal{F}^B_{q,t}(S)=\sum_{|\lambda|+|\rho|=n}m_{\lambda,\rho}(q,t)\, s_{\la}[\z+\bar{\z}] s_{\rho}[\z-\bar{\z}],
 \end{equation}
with
      $$m_{\lambda,\rho}(q,t):=\sum_{k,j\geq 0}  m_{\la,\rho}(k,j) q^k t^j,$$ 
where $m_{\la,\rho}(k,j)$ is the multiplicity of the irreducible representation of character $\charac^{\la,\rho}$ in the bihomogeneous component $S_{k,j}$. Using this Frobenious characteristic, it is easy to check that $h_n[{\bf z}+{\bf \bar{z}}]$ and $e_n[{\bf z}-{\bf \bar{z}}]$, correspond to the trivial and alternating representation, respectively.

For example, the character $\charac^{\QQ B_n}$ of the (left) regular representation of $B_n$ is readily seen to be 
    $$\charac^{\QQ B_n}(\mu,\nu)=\begin{cases}
    2^n n! & \text{if}\ (\mu,\nu)=(1^n,0), \\
     0 & \text{otherwise}.
\end{cases}$$
Hence, the corresponding Frobenius characteristic is
\begin{eqnarray}
     \mathcal{F}^B(\QQ B_n)&=& (2\, p_1(\z))^n \nonumber\\
         &=& (p_1(\z+\bar{\z}) +p_1(\z-\bar{\z}))^n\label{frob_regular}\\
         &=&\sum_{|\lambda|+|\rho|=n} {n\choose |\lambda|} f_\lambda f_\rho s_\lambda[\z+\bar{\z}] s_\rho[\z-\bar{\z}].\label{dim_irred}
\end{eqnarray}
Once again formula (\ref{dim_irred}) corresponds, now in the $B_n$ case, to the classical decomposition of the regular representation with multiplicities
of each irreducible character equal to its dimension.

\section{Frobenius characteristic of $\QQ[\x]$.}
It is not hard to show that the  graded Frobenius characteristic of the ring of polynomials $Q=\QQ[\x]$  is simply
\begin{equation}\label{bigraded_frob_Qx}
     \mathcal{F}^B_q(Q)=h_n\left[\frac{\z}{1-q}+\frac{\bar{\z}}{1+q}\right].
\end{equation}
To see this, we compute directly the value at $(\mu,\nu)$, of the graded character of $Q$,  in the basis of monomials. In broad terms, this computation goes as follows. Let $\beta$ in $B_n$ be of signed cycle type $(\mu,\nu)$, the only possible contribution to the trace of $\beta\cdot(-)$ has to come from monomials $\x^{\bf a}$ fixed up to sign by $\beta$. This forces the entries of $\bf a$ to be constant on cycles of $\beta$, and the sign is easily computed. It follows that
\begin{equation}
   \sum_{d\geq 0} \charac^{Q_d}(\mu,\nu)\, q^d=  
              \prod_{i=1}^{\ell(\mu)} \frac{1}{1-q^{\mu_i}}
              \prod_{j=1}^{\ell(\nu)} \frac{1}{1+q^{\nu_j}}.
\end{equation}
Hence from (\ref{h-to-p}) and (\ref{hn_sum}) we get
 \begin{eqnarray*}
     \mathcal{F}^B_q(Q)&=&\sum_{|\mu|+|\nu|=n} 
             \frac{1}{z_\mu} p_\mu\left[ \frac{\z}{1-q}\right]
              \frac{1}{z_\nu} p_\nu\left[ \frac{\bar{\z}}{1+q}\right]\\
          &=&\sum_{k=0}^n h_k\left[ \frac{\z}{1-q}\right] h_{n-k}\left[ \frac{\bar{\z}}{1+q}\right] \\
          &=&  h_n\left[\frac{\z}{1-q} +\frac{\bar{\z}}{1+q}\right],
\end{eqnarray*}  
Using (\ref{cauchy}), we can express this as
 \begin{equation}\label{frob_Q}
\mathcal{F}_q^B(Q) = \sum_{|\la|+|\rho|=n}s_{\la}\left[\frac{1}{1-q^2}\right]s_{\rho}\left[\frac{q}{1-q^2}\right] s_{\la}[\z+\bar{\z}] s_{\rho}[\z-\bar{\z}].
\end{equation}
In particular, we get back formula (\ref{hilbert_invariant})
 \begin{equation*}
     m_{(n),0}(q)=\prod_{i=1}^n \frac{1}{1-q^{2i}}
 \end{equation*}
in the special case of $B_n$.

\section{Frobenius characteristic of $\mathcal{H}$.}

We can now start our investigation of the $B_n$-module structure of the space of $B_n\times B_n$-harmonics 
     $$\mathcal{H}=\mathcal{H}_{B_n \times B_n}.$$ 
The first step is to compute its bigraded Frobenius characteristic, using the simply graded case as a stepping stone. From  (\ref{chev_harm}), (\ref{Hilbertinvariants}), (\ref{bigraded_frob_Qx}), and using the fact that invariant polynomials play the role of ``constants'' in the context of representation theory, we get
\begin{equation}\label{frobenius_harmx}
    \mathcal{F}^B_q(\mathcal{H}_{B_n})= h_n\left[\frac{\z}{1-q}+\frac{\bar{\z}}{1+q}\right] \prod_{i=1}^n (1-q^{2i}).
\end{equation}
In other words (in view of (\ref{frob_Q})),  the graded multiplicity of the irreducible components of type $(\la,\rho)$ is the polynomial\footnote{It is clearly a polynomial since $\mathcal{H}_{B_n}$ is finite dimensional.}
\begin{equation}
   m_{\lambda,\rho}(q)=  s_{\la}\left[\frac{1}{1-q^2}\right]s_{\rho}\left[\frac{q}{1-q^2}\right]
         \prod_{i=1}^n (1-q^{2i})
\end{equation}
It may be worth mentioning that there is a combinatorial interpretation  for $m_{\lambda,\rho}(q)$, in the form:
\begin{equation}
    m_{\lambda,\rho}(q)=\sum_{{\bf t}\in {\rm SYT}(\lambda,\rho)} q^{f({\bf t})}
\end{equation}
with ${\bf t}$ running over some natural set of ``standard Young tableaux'' of shape $(\lambda,\rho)$, and with $f({\bf t})$ a notion of ``major index'' for such tableaux. For more details see  (see \cite[Section 5]{brenti}).
With $n=3$, we get the following values for $m_{\lambda,\rho}(q)$:
\begin{center}
\begin{picture}(100,70)(0,0)
\put(-110,35){\framebox{$ \begin {array}{cccccccccc} 
(\lambda,\rho):&(3,0)&(21,0)&(111,0)&(2,1)&(11,1)\\\noalign{\medskip}
m_{\lambda,\rho}(q): &1&  {q}^{4}+q^2   & {q}^{6}&  {q}^{5}+{q}^{3}+q &  {q}^{7}+{q}^{5}+q^3 
        \\\noalign{\bigskip}
(\lambda,\rho):&(0,111)&(0,21)&(0,3)&(1,11)&(1,2)\\\noalign{\medskip}
m_{\lambda,\rho}(q): &{q}^{9}&{q}^{5}+q^7&{q}^{3}& {q}^{4}+{q}^{6}+q^8 &{q}^{2}+{q}^{4}+q^6   
\end{array}  $}}
\put(-110,37){\line(1,0){325}}
\end{picture}
\end{center}
\setlength{\unitlength}{4mm}
These specific values have been organized in order to make evident that
\begin{equation}\label{symetrie_flip}
    m_{\rho',\lambda'}(q)=q^{n^2} m_{\lambda,\rho}(1/q)
\end{equation}
This is a general fact to which we will come back later.
Now, it is easy to observe that the character of a tensorial product is the product of the characters of each terms. Thus, for two $B_n$-modules $\mathcal{V}$ and $\mathcal{W}$, we have
\begin{equation}\label{frob_tensor}
    \mathcal{F}^B(\mathcal{V}\otimes \mathcal{W})=\mathcal{F}^B(\mathcal{V}) * \mathcal{F}^B(\mathcal{W}),
\end{equation}
where $``*"$ stands for {\em internal product} of symmetric functions. Recall that this is the bilinear product on $\Lambda(\z)\otimes \Lambda(\bar{\z})$ such that
\begin{equation}\label{internal}
     p_{\mu,\nu}(\z,\bar{\z})* p_{\la,\rho}(\z,\bar{\z})
     =\begin{cases}
      z_\mu\, z_\nu\,p_{\mu,\nu}(\z,\bar{\z}) 
             & \text{if\ }(\mu,\nu)=(\lambda,\rho).  \\
      0 & \text{otherwise}.
\end{cases}
\end{equation}
This last definition immediately implies that     
    $$f(\z) *g(\z) =0,\qquad {\rm and}\qquad f(\bar{\z}) *g(\bar{\z}) =0,$$
whenever $f$ and $g$ are of different homogeneous degree. Moreover (\ref{internal}) implies 
\begin{equation}\label{prod-omega}
    \Omega[a \z] \Omega[b {\bf \bar z}]* \Omega[c\z] \Omega[bd {\bf \bar z}]= \Omega[ac \z] \Omega[bd {\bf  \bar z}].
\end{equation}
In particular, $\Omega[\z+\bar{\z}]$ is the neutral element for the internal product (\ref{internal}).
We can then apply (\ref{frobenius_harmx}) to easily compute the bigraded Frobenius characteristic of $\mathcal{H}$. By (\ref{d-H}) and (\ref{frob_tensor}) we get $$\mathcal{F}_{q,t}^B(\mathcal{H})=\mathcal{F}_{q}^B(\mathcal{H}_{B_n})* 
                    \mathcal{F}_{t}^B(\mathcal{H}_{B_n}).$$
Then by using (\ref{omega}) and (\ref{prod-omega}) we obtain
\begin{eqnarray*}
   \Omega\left[\frac{\z}{1-q}+\frac{\bar{\z}}{1+q}\right]&*&
       \Omega\left[\frac{\z}{1-t}+\frac{\bar{\z}}{1+t}\right]  =\\
       &&\quad \Omega\left[\frac{\z}{(1-q)(1-t)}+\frac{\bar{\z}}{(1+q)(1+t)}\right].
 \end{eqnarray*}
This readily implies 
\begin{eqnarray}
     \mathcal{F}_{q,t}^B(\mathcal{H})&=&   
        h_n\left[\frac{\z}{(1-q)(1-t)}+\frac{\bar{\z}}{(1+q)(1+t)}\right]\times
           \prod_{i=1}^n(1-q^{2i})(1-t^{2i}) \label{frob_R} \nonumber\\
         &=& \sum_{|\lambda|+|\rho|=n} s_{\lambda}\left[\frac{1+qt}{(1-q^2)(1-t^2)}\right]   s_{\rho}\left[\frac{t+q}{(1-q^2)(1-t^2)}\right] s_{\lambda}[\z + \bar{\z}] s_{\rho}[\z - \bar{\z}] \nonumber \\
         & & \quad  \times
           \prod_{i=1}^n(1-q^{2i})(1-t^{2i}). 
 \end{eqnarray}             
It follows, in particular,  that the corresponding graded multiplicity of  the trivial representation is
\begin{equation}\label{mult_trivial}
     H_{q,t}(\mathcal{H}^{B_n})=h_n\left[\frac{1+q\,t}{(1-q^2)(1-t^2)}\right] \prod_{i=1}^n (1-q^{2i})\prod_{i=1}^n (1-t^{2i})
\end{equation}
where we use the notation $H_{q,t}$ for the {\em bigraded Hilbert series}.
Expanding $h_n$ in power sum, it is easy to show that (\ref{frob_R}) specializes, when $t$ tends to $1$,  to
   $$\lim_{t\rightarrow 1} \mathcal{F}_{q,t}^B(\mathcal{H})=
          (p_1(\z+\bar{\z})+p_1(\z-\bar{\z}))^n \prod_{j=1}^n \frac{1-q^{2j}}{1-q},$$
 which is clearly a multiple of the character (\ref{frob_regular}) of the regular representation of  $B_n$. By symmetry, a similar expression in $t$ results when we rather let $q$ tend to $1$ in (\ref{frob_R}). An immediate corollary, is that each irreducible representation $\mathcal{V}^{\lambda,\rho}$, of $B_n$, appears in $\mathcal{H}$ with multiplicity (see (\ref{dim_irred}))  equal to
     $$2^n\, n!\, {n\choose |\lambda|}f_{\lambda}f_{\rho}.$$

\section{The trivial component of $\mathcal{C}$}\label{triv-section}
In view of (\ref{mult_trivial}), and the discussion that follows, the dimension of trivial isotypic component of $\mathcal{H}$ as well as that of $\mathcal{C}$ is equal to the order of $B_n$. It is thus natural to expect the existence of a basis of $\mathcal{C}^{B_n}$ indexed naturally by elements of $B_n$. We will describe such a basis in Section \ref{base}. To this end, we need to introduce some notation.

Let both ${\bf a}$ and ${\bf b}$ be length $n$ vectors of nonnegative integers. The pair $({\bf a},{\bf b})$ is said to be a {\em bipartite vector}, and we will often use the {\em two-line notation}
\begin{equation*}
     ({\bf a},{\bf b})= \begin{pmatrix} a_1 & a_2 & \ldots & a_n\\ 
                       b_1 & b_2 & \ldots & b_n \end{pmatrix}.
\end{equation*}
Along the same lines, we say that a bipartite vector $({\bf a},{\bf b})$ is a {\em bipartite partition} if
its {\em parts}, $(a_i,b_i)$, are ordered in increasing {\em lexicographic order} (see also \cite{garsia_gessel}, \cite{ges}). This is to say the order,   ``$\prec_{lex}$'', for which
 $$(a,b)\prec_{lex} (a',b') \iff \begin{cases}
      a < a'  & \text{or}, \\
      a=a' & \text{and}\quad b<b'.
\end{cases}$$
Notice that all are partitions are of length $n$, and the we are allowing empty parts $(0,0)$, all situated at the beginning.
For any bipartite vector $({\bf a},{\bf b})$, we denote $({\bf a},{\bf b})^\circlearrowleft$ the bipartite partition obtained by reordering the pairs $(a_i,b_i)$ in increasing lex-order.

Now, setting
\begin{equation*}
    x_1<y_1<x_2<y_2<\ldots <x_n<y_n,
\end{equation*}
 as the underlying order for the variables,  we can make a Gr\"obner basis argument  very similar to the one used for (\ref{basis_coinvariant}), to deduce that the $(2^n n!)^2$ dimensional space of coinvariants $\mathcal{C}=R_{B_n\times B_n}$ affords the basis
  $$\{{\bf x^a}{\bf y^b}+\mathcal{I}_{B_n\times B_n}\ | \ a_i<2i \  {\rm and} \   b_j<2j\}.$$
Thus $\mathcal{H}=\mathcal{H}_{B_n\times B_n}$, the corresponding space of harmonics, affords the linear basis
  $$\{\partial{\bf x^a}\partial{\bf y^b} \Delta(\x)\Delta(y)\ | \ a_i<2i \  {\rm and} \   b_j<2j\},$$
 since the Jacobian determinant of $B_n\times B_n$ is clearly $\Delta(\x)\Delta(\y)$ (see (\ref{jacobian}).  
Without surprise, it develops that a linear basis for the set of diagonally invariant polynomials $R^{B_n}$ is naturally indexed by  {\em even} bipartite partitions. These are bipartite partitions 
\begin{equation*}
   D=({\bf a},{\bf b})
\end{equation*}
such that each $a_i + b_i $ is even. In the sequel we will use the term {\em $e$-diagram} (rather then that of ``even bipartite partition'') for reasons that will become clear in the next section.
The announced homogeneous basis for $R^{B_n}$ is simply given by the set 
  $$\{M{({\bf a},{\bf b})}\ |\ ({\bf a},{\bf b})\ {\rm is \; an} \;   \hbox{\rm {\em e}-diagram},\          
                \ell({\bf a},{\bf b})= n\},$$
of {\em monomial} diagonal invariants, defined as:
   $$M{({\bf a,b})}:=\sum\{ \sigma\cdot \x^{\bf a}\y^{\bf b} \ |\ {\sigma\in S_n}\},$$
where we sum over the set of (hence distinct) monomials obtained by permuting the variables. For example, with $n=3$,
 $$ M{\left( \begin{smallmatrix} 0 &0 & 2 \\
                        0 & 4 & 2  \end{smallmatrix}\right)}
    =y_2^4x_3^2y_3^2+x_2^2y_2^2y_3^4 +y_1^4x_3^2y_3^2 +y_1^4x_2^2y_2^2 +x_1^2y_1^2y_3^4 + x_1^2y_1^2y_2^4.$$
Observe that the leading monomial of $M{({\bf a,b})}$ is ${\bf x^{a}y^{b}}$. Clearly, each $M{({\bf a,b})}$ is a bihomogeneous polynomial of bidegree $(|{\bf a}|,|{\bf b}|)$. We will also say that
\begin{equation*}
    |D|:=(|{\bf a}|,|{\bf b}|),
\end{equation*}
is the {\em weight} of the $e$-diagram $D$.

To construct a basis of $\mathcal{C}^{B_n}$, naturally indexed by the group elements, we associate to each $\be$ in $B_n$, a special $e$-diagram $D_{\be}$ as follows.  Closely imitating  (\ref{fmaj_i}), (\ref{varepsilon}) and (\ref{desc_i}),   let us set for each $i$,
\begin{eqnarray}\label{def-g_i}
    g_i(\be)&:=&2 \de_i(\be)+\eta_i(\be),\  {\rm with}\\
   \eta_i(\beta)&:=& \begin{cases}
     1 & \textrm{if }\be(i) > 0,\ {\rm and}\nonumber \\
      0 & \text{otherwise},\nonumber
\end{cases}\\
    \de_i(\beta)&:=&\# \{j \in \Des(\be)\ |\ j <i \}.\nonumber
\end{eqnarray}
Clearly,
\begin{equation}\label{cresce}
g_i(\be)\leq g_{i+1}(\be).
\end{equation}
We then define the $e$-diagram 
\begin{equation}\label{D_beta}
    D_{\be}:=({\bf g}(\be),\widetilde{{\bf g}}(\be)),
 \end{equation}
 with
\begin{equation}\label{defediagram}
\begin{array}{rcl}
     {\bf g}(\be)&:=&(g_1(\be),g_2(\be),\ldots, g_n(\be)),\\
     \widetilde{{\bf g}}(\be)&:=& ( \widetilde{g}_1(\be), \widetilde{g}_2(\be),\ldots, \widetilde{g}_n(\be)), \quad {\rm with}\quad
         \widetilde{g}_i(\be):=g_{\sigma(i)}(\be^{-1}),
\end{array}\end{equation}
and where $\sigma(i)$ equal to the absolute value of $\beta(i)$.
For example, with $n=2$, we get
\begin{equation}\label{diaga_B2}
\qquad     \begin{matrix}
          D_{12}=\left( \begin{smallmatrix} 1& 1\\ 1& 1 \end{smallmatrix}\right), &
          D_{\bar{1}2}=\left( \begin{smallmatrix} 0& 1\\ 0& 1 \end{smallmatrix}\right), &
          D_{1\bar{2}}=\left( \begin{smallmatrix} 1& 2\\ 1& 2 \end{smallmatrix}\right), &
          D_{\bar{1}\bar{2}}=\left( \begin{smallmatrix} 0& 0\\ 0& 0 \end{smallmatrix}\right),  \\ \\
         D_{21}=\left( \begin{smallmatrix} 1& 3\\ 3& 1 \end{smallmatrix}\right), &
          D_{\bar{2}1}=\left( \begin{smallmatrix} 0& 1\\ 2& 1 \end{smallmatrix}\right), &
          D_{2\bar{1}}=\left( \begin{smallmatrix} 1& 2\\ 1& 0 \end{smallmatrix}\right), &
          D_{\bar{2}\bar{1}}=\left( \begin{smallmatrix} 0& 2\\ 2& 0 \end{smallmatrix}\right),  
     \end{matrix}
\end{equation}
 We will show in Section \ref{base} that every element in $R^{B_n}$ can be written in terms of the set 
\begin{equation}\label{base_M}
  \mathcal{M}_n:=\{M_\beta\ |\ \beta\in B_n\},
\end{equation}
using the simpler notation $M_\beta$, instead of $M(D_\beta)$, for the monomial diagonal invariants associated to these special diagrams. For example, following (\ref{diaga_B2}) the 8 elements of $\mathcal{M}_2$ are
$$\begin{array}{ll}
M_{12}=x_1y_1\, x_2y_2,& M_{\bar{2}\bar{1}}=y_1^2\, x_2^2 + x_1^2\,y_2^2,\\ \\
M_{\bar{1}2} =x_1y_1+x_2y_2,&M_{1\bar{2}}=x_1y_1\, x_2^2y_2^2 + x_1^2y_1^2\, x_2y_2, \\ \\
 M_{2\bar{1}}=x_1y_1 \, x_2^2+ x_1^2\, x_2y_2,&M_{\bar{2}1}=y_1^2\, x_2y_2 + x_1y_1\, y_2^2,\\ \\
 M_{\bar{1}\bar{2}}=1,&M_{21}=x_1y_1^3\, x_2^3y_2 + x_1^3y_1\, x_2y_2^3.
 \end{array}$$
More precisely, we will see later (see Section \ref{base}) that, for every $e$-diagram 
$({\bf c},{\bf d})$,  there are unique  $u_{\be}(\x,\y)$'s in $ \QQ[\x]^{B_n} \otimes \QQ[\y]^{B_n}$ such that
\begin{equation}\label{decomposition_unique}
   M{({\bf c},{\bf d})}=\sum_{\be \in B_n} u_{\be}(\x,\y) \; M_\beta.
 \end{equation}
 Recall that we are writing $M_\beta$ instead of $M{({\bf g}(\be), \widetilde{{\bf g}}(\be))}$.
This immediately implies that the Hilbert series of $\mathcal{M}_n$ equals $H_{q,t}(\mathcal{C}^{B_n})$, for which we have obtained expression (\ref{mult_trivial}). 
Now, consider the involution on $B_n$, $\be \mapsto \be ^{\circ}$, such that
\[\be^{\circ}:=- w_{0} \be\, w_0,\]
where $w_0:=n\cdots 2\,1$, and the minus stands a global sign change. Evidently, we have 
\begin{eqnarray*} 
 \de_i(\be)&=&d_{n+1-i}(\be^{\circ})\\
 \eta_i(\be)&=&\varepsilon_{n+1-i}(\be^{\circ});
\end{eqnarray*}
and using the relevant definitions (\ref{defediagram}) and (\ref{fmaj}), it follows that
  $$|{\bf g}(\be)|=\fmaj(\be^{{\circ}}) \quad  {\rm and} \quad  |\widetilde{{\bf g}}(\be)|=\fmaj(\be^{-1})^{{\circ}}.$$
We can thus conclude, modulo a proof of (\ref{decomposition_unique}), that we have the bigraded Hilbert series formula
\begin{equation}\label{hilbert_compact}
H_{q,t}(\mathcal{C}^{B_n})=\sum_{\be \in B_n}q^{\fmaj(\be)}t^{\fmaj(\be^{-1})}.
\end{equation}
From this and (\ref{the_iso_inv}), it follows the graded Hilbert series of the diagonally invariant polynomials is given by  
\begin{equation}
H_{q,t}(R^{B_n})= h_n\left[\frac{1+q\,t}{(1-q^2)(1-t^2)}\right]=\sum_{\be \in B_n}q^{\fmaj(\be)}t^{\fmaj(\be^{-1})}
 \prod_{i=1}^n \frac{1}{1-q^{2i}}\prod_{i=1}^n \frac{1}{1-t^{2i}}\label{genfunction}.
\end{equation}
For different approaches to this computation we refer the reader to \cite{adin} and \cite{bc}. 

\section{Combinatorics of $e$-diagrams}\label{combin}
We are now going to unfold a combinatorial approach to the concepts introduced in the previous sections. Central to our discussion is a natural classification of $e$-diagrams in term of elements of $B_n$, envisioning them as (multi) subsets of the {\em even-chessboard plane}. This is simply the $(\NN\times \NN)^0$ subset of the combinatorial plane defined as
   $$(\NN\times \NN)^0:=\{(a,b) \in \NN \times \NN\ |\quad a+b\equiv 0\ ({\rm mod}\ 2) \}.$$
 Its elements are called {\em cells}.
Clearly, $e$-diagrams  correspond to $n$-cell multisubsets in $(\NN\times \NN)^0$. Here, the word  ``multisubset'' underlines that we are allowing cells to have multiplicities.

\begin{figure}
\begin{center}\begin{picture}(11,11)(0,0)
\multiput(1,0)(2,0){5}{\carre}
\multiput(0,1)(2,0){5}{\carre}
\multiput(1,2)(2,0){5}{\carre}
\multiput(0,3)(2,0){5}{\carre}
\multiput(1,4)(2,0){5}{\carre}
\multiput(0,5)(2,0){5}{\carre}
\multiput(1,6)(2,0){5}{\carre}
\multiput(0,7)(2,0){5}{\carre}
\multiput(1,8)(2,0){5}{\carre}
\multiput(0,9)(2,0){5}{\carre}
\put(0.5,0.03){\jaune{\circle*{.9}}\hskip-10pt\rouge{\circle{.9}}}\put(0.25,-.25){$\rouge{\mathbf{2}}$}
\put(1.5,5.03){\jaune{\circle*{.9}}\hskip-10pt\rouge{\circle{.9}}}\put(1.25,4.75){$\rouge{\mathbf{1}}$}
\put(2.5,6.03){\jaune{\circle*{.9}}\hskip-10pt\rouge{\circle{.9}}}\put(2.25,5.75){$\rouge{\mathbf{2}}$}
\put(6.5,4.03){\jaune{\circle*{.9}}\hskip-10pt\rouge{\circle{.9}}}\put(6.25,3.75){$\rouge{\mathbf{1}}$}
\put(9.5,5.03){\jaune{\circle*{.9}}\hskip-10pt\rouge{\circle{.9}}}\put(9.25,4.75){$\rouge{\mathbf{1}}$}
\put(8.5,0.03){\jaune{\circle*{.9}}\hskip-10pt\rouge{\circle{.9}}}\put(8.25,-.25){$\rouge{\mathbf{1}}$}
\put(9.5,9.03){\jaune{\circle*{.9}}\hskip-10pt\rouge{\circle{.9}}}\put(9.25,8.75){$\rouge{\mathbf{1}}$}
\thicklines
\multiput(0,-0.5)(0,2){6}{\line(1,0){10}}
\multiput(0,-0.5)(2,0){6}{\line(0,1){10}}
\thinlines
\multiput(0,0.5)(0,2){5}{\line(1,0){10}}
\multiput(1,-0.5)(2,0){5}{\line(0,1){10}}
\end{picture}
\caption{A chessboard diagram.}\label{diagramma}
\end{center}
\end{figure}
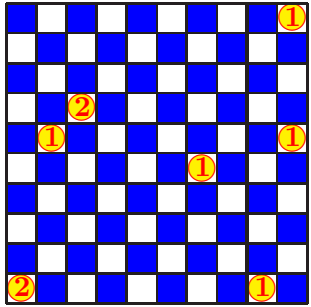

\noindent Figure \ref{diagramma} gives the graphical representation of the $e$-diagram
\begin{equation}\label{un_diag}
    \begin{pmatrix} 0 & 0 & 1  & 2 & 2 & 6 & 8 & 9 & 9 \\
                             0 & 0 & 5  & 6 & 6 & 4 & 0 & 5 & 9 \end{pmatrix}
 \end{equation}
 that is implicit in this ``geometrical'' outlook.  Integers in the cells give the associated multiplicity.

\noindent The set of  {\em descents}, $\Des(D)$, and the set of {\em sign changes}, $\Sx(D)$,  of an $e$-diagram 
     $$D=\begin{pmatrix} a_1 & a_2 & \ldots & a_n\\ 
                       b_1 & b_2 & \ldots & b_n \end{pmatrix}$$
are defined to be
\begin{eqnarray*}
     \Des(D)&:=&\{k \in [n-1]\ |\ b_k>b_{k+1} \  {\rm and} \  a_k\equiv a_{k+1}\ ({\rm mod}\ 2)\},\\
    \Sx(D)&:=&\{k \in [0,n-1]\ |\  a_k \not  \equiv a_{k+1}\ ({\rm mod}\ 2)\},
 \end{eqnarray*}
where by convention  $a_0\equiv 1\ ({\rm mod}\ 2)$. In other words, $0$ is in $\Sx(D)$ if and only if $a_1$ is odd. Note that these two set are disjoint, i.e.,
  $$\Des(D) \cap \Sx(D) =\emptyset.$$
For the $e$-diagram of Figure \ref{diagramma}, we have
  $$\Des(D)=\{5,6\}\qquad {\rm and}\qquad \Sx(D)=\{2,3,7\}.$$
We further set
\begin{eqnarray*}
   g_i(D)&:=&2 \de_i(D)+\sx_i(D), \quad {\rm with}\\
   \de_i(D)&:=& \# \{k \in \Des(D)\ |\ k <i\ \},\quad {\rm and}  \\
    \sx_i(D)&:=& \# \{k \in \Sx(D)\ |\ k  <i\ \} .
 \end{eqnarray*}
Then, for any cell $(a_i,b_i)$ in $D$, we must have
\begin{equation}\label{piugrande}
a_i\geq g_i(D).
\end{equation}
We now associate to each $e$-diagram $D$ a signed permutation $\be(D)\in B_n$.
In order to define this $\be(D)$ we suppose that the cells of
     $$D=\begin{pmatrix} a_1 & a_2 & \ldots & a_n\\ 
                       b_1 & b_2 & \ldots & b_n \end{pmatrix}$$
have been ordered in increasing {\em opposite lexicographical order}. This is say the order, denoted ``$\prec_{op}$'', such that
 $$(a,b)\prec_{op} (a',b') \iff\begin{cases}
      b < b' \quad \text{or}, \\
      b = b'\quad {\rm and}\quad a<a'.
     \end{cases} $$
We will call this the {\em labelling order} for cells of the diagram. Our intent here is that a cell $(a,b)$ of $D$ be labeled $i$, if it sits in $D$ in the $i^{\rm th}$ position  with respect to this labeling order. There is clearly a unique permutation $\sigma$ such that
\[(a_{\s(1)},b_{\s(1)})\preceq_{lex}(a_{\s(2)},b_{\s(2)})\preceq_{lex} \ldots \preceq_{lex}(a_{\s(n)},b_{\s(n)}),\]
with $\s(i)<\s(j)$, whenever $i<j$ and $(a_i,b_i)=(a_j,b_j)$. We then introduce the  {\em  classifying signed permutation}, $\beta=\be(D)$, of $D$, setting
        $$\beta(i):=(-1)^{a_i+1} \sigma(i).$$
In this context, we often say that $\prec_{lex}$ is the {\em reading order} for the cells of $D$. 
For the diagram of (\ref{un_diag}), represented in Figure \ref{diagramma}, the labeling order is   
   $$\begin{pmatrix} 0 & 0 & 8 & 6 & 1 & 9  & 2 & 2  & 9 \\
                                0 & 0 & 0 & 4 & 5 & 5  & 6 & 6  & 9 \end{pmatrix}.$$ 
The corresponding classifying signed permutation is readily seen to be $\bar{1}\bar{2}5 \bar{7}\bar{8} \bar{4}\bar{3} 6 9 .$

A simple inductive argument shows that for all $i$'s,
\begin{equation}\label{hi}
    g_i(D)=g_i(\be(D)),
\end{equation}
justifying our use of the same notation in both cases.
Now, if $D=({\bf a},{\bf b})$, let us set
\begin{equation}\label{inversdiagram}
    D^{*}:=({\bf b},{\bf a})^\circlearrowleft.
\end{equation}
Recall that $\circlearrowleft$  indicates that we are passing to the associated bipartite partition. This is to say that we are reordering the cells in increasing lex-order.
It is easy to check that
  $$\be(D^*)=\be^{-1},$$
 if $\beta=\beta(D)$.
Thus, we naturally call $D^*$ the {\em inverse e-diagram} of $D$.
It follows from (\ref{hi}) that 
\begin{equation}\label{hinversa}
g_i(D^*)=g_i(\be^{-1}).
\end{equation}
Hence
\begin{equation}\label{piu_vert}
      b_i\geq g_{\sigma(i)}(D^*),
\end{equation}
where, as before, $\sigma(i)$ is the absolute value of $\beta(i)$.
All this suggests that, in a manner similar to (\ref{defediagram}), we use the notation
$$\begin{array}{rcl}
     {\bf g}(D)&:=&(g_1(D),g_2(D),\ldots, g_n(D)),\\
     \widetilde{{\bf g}}(D)&:=&( \widetilde{g}_1(D), \widetilde{g}_2(D),\ldots, \widetilde{g}_n(D)), \quad {\rm with}\quad
         \widetilde{g}_i(D):=g_{\sigma(i)}(D^*),
\end{array}$$
where, once again,  $\sigma(i)$ is  the absolute value of $\beta(i)$. We have thus associated to each diagram, $D$, a new diagram 
\begin{equation}\label{compactifie}
     \overline{D}:=({\bf g}(D),\widetilde{{\bf g}}(D)),
 \end{equation}
which is going to be called\footnote{For reasons that will be made clear in the next section.} the {\em compactification} of $D$.
The compact diagram associated to (\ref{un_diag}) is
\begin{equation}\label{un_diag_comp}
    \begin{pmatrix} 0 & 0 & 1  & 2 & 2 & 4 & 6 & 7 & 7 \\
                             0 & 0 & 3  & 4 & 4 & 2 & 0 & 3 & 5 \end{pmatrix}.
 \end{equation}
We can check, and this will be made clear in the next Section \ref{sec_comp}, that
\begin{prop} For all  $e$-diagram $D$, we have
\begin{equation}\label{dbeta}
       \be(D)=\be(\overline{D});
 \end{equation}
 and 
\begin{equation}\label{equazioni}
      \overline{D}=D_\beta,
\end{equation}
with $\be:=\be(D)$.
\end{prop}
We introduce an equivalence relation on the set of the $e$-diagrams saying that $D$ and $\tilde{D}$ are equivalent if and only if they have the same classifying signed permutation. In symbols,
    $$D \simeq \widetilde{D} \iff \be(D)=\be(\widetilde{D}).$$
In view of (\ref{cresce}) the cells of $D_{\beta}$ are in reading order. Moreover the label of $(g_i(\beta),\widetilde{g}_i(\beta))$ is $\beta(i)$, hence $\beta(D_{\beta})=\beta$.
For the moment, the main properties of this equivalence relation is the following.
\begin{thm} For any $e$-diagram, $D$, we have
    $$\be(D)=\be \Longleftrightarrow D \simeq D_{\be}.$$
Moreover, $D_{\be}$ is minimal  in the equivalence class of $D$, in that the matrix
   $$  \widetilde{D}-D_{\be}$$
has all entries nonnegative, for all $\widetilde{D} \simeq D$.
\end{thm}
The last part of this theorem is one of the reasons why we say that $e$-diagrams of the form $D_\beta$ are {\em compact}. The next section will make this notion even more precise.

\begin{figure}
\begin{center}
\begin{picture}(11,11)(0,0)
\multiput(1,0)(2,0){5}{\carre}
\multiput(0,1)(2,0){5}{\carre}
\multiput(1,2)(2,0){5}{\carre}
\multiput(0,3)(2,0){5}{\carre}
\multiput(1,4)(2,0){5}{\carre}
\multiput(0,5)(2,0){5}{\carre}
\multiput(1,6)(2,0){5}{\carre}
\multiput(0,7)(2,0){5}{\carre}
\multiput(1,8)(2,0){5}{\carre}
\multiput(0,9)(2,0){5}{\carre}
\multiput(1,5)(2,0){1}{\jcarre}\multiput(0,5)(2,0){2}{\vcarre}
\multiput(0,4)(2,0){5}{\jcarre}\multiput(1,4)(2,0){5}{\vcarre}
\multiput(5,3)(2,0){3}{\jcarre}\multiput(4,3)(2,0){3}{\vcarre}
\put(3.5,5){\circle*{.8}}
\put(3.5,3){\circle*{.8}}
\thicklines
\put(3.5,5){\linethickness{2pt}\rose{\vector(0,-1){1.5}}}
\multiput(0,-0.5)(0,2){6}{\line(1,0){10}}
\multiput(0,-0.5)(2,0){6}{\line(0,1){10}}
\thinlines
\multiput(0,0.5)(0,2){5}{\line(1,0){10}}
\multiput(1,-0.5)(2,0){5}{\line(0,1){10}}
\end{picture}\qquad
\begin{picture}(11,11)(0,0)
\multiput(1,0)(2,0){5}{\carre}
\multiput(0,1)(2,0){5}{\carre}
\multiput(1,2)(2,0){5}{\carre}
\multiput(0,3)(2,0){5}{\carre}
\multiput(1,4)(2,0){5}{\carre}
\multiput(0,5)(2,0){5}{\carre}
\multiput(1,6)(2,0){5}{\carre}
\multiput(0,7)(2,0){5}{\carre}
\multiput(1,8)(2,0){5}{\carre}
\multiput(0,9)(2,0){5}{\carre}
\multiput(5,1)(0,2){1}{\jcarre}\multiput(5,0)(0,2){2}{\vcarre}
\multiput(4,0)(0,2){5}{\jcarre}\multiput(4,1)(0,2){5}{\vcarre}
\multiput(3,5)(0,2){3}{\jcarre}\multiput(3,4)(0,2){3}{\vcarre}
\put(5.5,3){\circle*{.8}}
\put(3.5,3){\circle*{.8}}
\thicklines
\put(5.5,3){\linethickness{2pt}\rose{\vector(-1,0){1.5}}}
\multiput(0,-0.5)(0,2){6}{\line(1,0){10}}
\multiput(0,-0.5)(2,0){6}{\line(0,1){10}}
\thinlines
\multiput(0,0.5)(0,2){5}{\line(1,0){10}}
\multiput(1,-0.5)(2,0){5}{\line(0,1){10}}
\end{picture}
\caption{Constraints on compacting moves.}\label{allowable}
\end{center}
\end{figure}
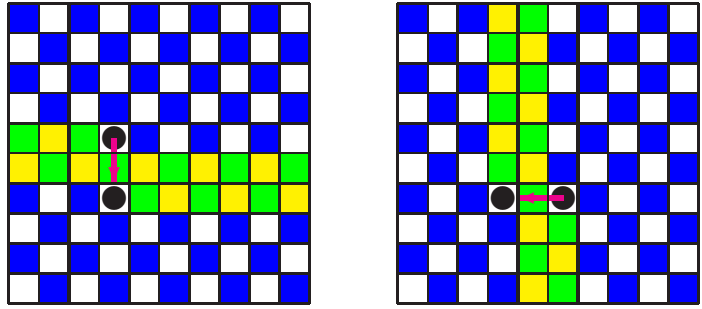

\section{Compactification of $e$-diagrams}\label{sec_comp}
As we will currently see, compact $e$-diagrams $D_{\be}$ can be characterized in a very simple ``geometrical'' manner. To obtain this characterization, we introduce the notion of ``compacting moves'' for  a diagram. These moves are designed to preserve the underlying classifying signed permutation. Moreover, they produce smaller diagram for the partial order
     $$D\leq \widetilde{D}.$$
This last statement means that the matrix $ \widetilde{D}-D$ only has nonnegative entries.
 
Now, let $D$ be an $e$-diagram in which we select some cell $c=(a,b)$, with $a\geq 2$. A {\em left move} of the cell $c=(a,b)$ in $D$ is defined to be
\begin{equation}\label{left_move}
     \triangleleft_c(D):=(D\setminus \{c\}) \cup \{(a-2,b)\}.
\end{equation}
Since it is understood here that cells are counted with multiplicities, the {\em set difference} and {\em union}, in the right hand side of (\ref{left_move}), are to be understood as multiset operations. Thus,  the result corresponds to decreasing by $1$ the multiplicity of $c$ in $D$, and increasing that of $(a-2,b)$ by $1$.
Now, if $c=(a,b)$, with $b\geq 2$, a  {\em down move} of the cell $c=(a,b)$ in $D$ is defined to be
\begin{equation}\label{down_move}
     \triangledown_c(D):=(D\setminus \{c\}) \cup \{(a,b-2)\}.
\end{equation}
{\em Compacting moves}, on a diagram $D$, are either left moves or down moves with some constraint on the choice of $c$ as described below.
A left move for $c$ is allowed as a compacting move if and only if the set 
\begin{equation}\label{vert}
\Ver(c,D):=\{c' \in D\ |\ \triangleleft(c)\prec_{lex} c' \prec_{lex} \triangledown(c)\}
\end{equation}
is empty.
Analogously, a down move for $c$ is allowed if and only if  the set 
\begin{equation}\label{horiz}
    \Hor(c,D):=\Ver(c^*,D^*)^*
 \end{equation}
is empty. The ``constraint''  intervals $\Ver$ and $\Hor$ are illustrated in Figure \ref{allowable}.

Observe that compacting moves do not change the sign of the cell that is moved. But more importantly, they do not change the classifying signed permutation of the $e$-diagram, namely,
 \begin{equation}\label{preserve_perm}
     D \simeq \triangleleft_c(D)  \simeq \triangledown_c(D),
 \end{equation}
whenever $c$ is so corresponds to a compacting move. The final crucial property of compacting moves is that
 \begin{equation}\label{reduce_weight}
     \triangleleft_c(D)<D, \qquad {\rm and} \qquad \triangledown_c(D)<D.
\end{equation}
Diagram for which no compacting move are possible will be called {\em compact}. From (\ref{piugrande}) and (\ref{piu_vert}), it follows that the set of compact $e$-diagrams (see (\ref{D_beta})) is exactly
     $$\{\ D_{\be}\ |\ \beta\in B_n\ \}.$$
Observe that, starting with any given $e$-diagram $D$, if one keeps applying compacting moves (in whatever order) until no such move is possible, then the final result will always be the compact diagram $\overline{D}$.
See Figure \ref{B2} for all eight compact $e$-diagrams associated to elements of  $B_2$.
In our upcoming discussion, it will be helpful to organize compacting moves in groups, called {\em big compacting moves}. This is natural in light of the following observation. Whenever a left compacting move is possible for a cell $c=(a_i,b_i)$, then all cells that are larger then $c$, in reading order, will eventually be left moved in the compacting process. We may as well achieve all this in one step:
 $$\begin{pmatrix} a_1 &\ldots & a_i & \ldots & a_n\\ 
               b_1 & \ldots & b_i & \ldots & b_n \end{pmatrix} \rightsquigarrow
       \begin{pmatrix} a_1 &\ldots & a_i-2  & \ldots & a_n-2\\ 
               b_1  &\ldots & b_i  & \ldots & b_n \end{pmatrix}.    $$

\setlength{\unitlength}{4mm}
\def\carre{\bleu{\linethickness{\unitlength}\line(1,0){1}}}
\begin{figure}
$$\begin{array}{cccc}
\begin{picture}(4.5,4.5)(0,-.4)
\put(1,0){\carre}\put(3,0){\carre}
\put(0,1){\carre}\put(2,1){\carre}
\put(1,2){\carre}\put(3,2){\carre}
\put(0,3){\carre}\put(2,3){\carre}
\put(1.5,1.03){\jaune{\circle*{.9}}\hskip-10pt\rouge{\circle{.9}}}
\put(1.25,.75){$\rouge{\mathbf{2}}$}
\thicklines
\multiput(0,-0.5)(0,2){3}{\line(1,0){4}}
\multiput(0,-0.5)(2,0){3}{\line(0,1){4}}
\thinlines
\multiput(0,0.5)(0,2){2}{\line(1,0){4}}
\multiput(1,-0.5)(2,0){2}{\line(0,1){4}}
\end{picture} &
\begin{picture}(4.5,4.5)(0,-0.4)
\put(1,0){\carre}\put(3,0){\carre}
\put(0,1){\carre}\put(2,1){\carre}
\put(1,2){\carre}\put(3,2){\carre}
\put(0,3){\carre}\put(2,3){\carre}
\put(1.5,1.03){\jaune{\circle*{.9}}\hskip-10pt\rouge{\circle{.9}}}
\put(1.25,.75){$\rouge{\mathbf{1}}$}
\put(0.5,0.03){\jaune{\circle*{.9}}\hskip-10pt\rouge{\circle{.9}}}
\put(0.25,-0.25){$\rouge{\mathbf{1}}$}
\thicklines
\multiput(0,-0.5)(0,2){3}{\line(1,0){4}}
\multiput(0,-0.5)(2,0){3}{\line(0,1){4}}
\thinlines
\multiput(0,0.5)(0,2){2}{\line(1,0){4}}
\multiput(1,-0.5)(2,0){2}{\line(0,1){4}}
\end{picture} &
\begin{picture}(4.5,4.5)(0,-0.4)
\put(1,0){\carre}\put(3,0){\carre}
\put(0,1){\carre}\put(2,1){\carre}
\put(1,2){\carre}\put(3,2){\carre}
\put(0,3){\carre}\put(2,3){\carre}
\put(2.5,2.03){\jaune{\circle*{.9}}\hskip-10pt\rouge{\circle{.9}}}
\put(2.25,1.75){$\rouge{\mathbf{1}}$}
\put(1.5,1.03){\jaune{\circle*{.9}}\hskip-10pt\rouge{\circle{.9}}}
\put(1.25,.75){$\rouge{\mathbf{1}}$}
\thicklines
\multiput(0,-0.5)(0,2){3}{\line(1,0){4}}
\multiput(0,-0.5)(2,0){3}{\line(0,1){4}}
\thinlines
\multiput(0,0.5)(0,2){2}{\line(1,0){4}}
\multiput(1,-0.5)(2,0){2}{\line(0,1){4}}
\end{picture} &
\begin{picture}(4.5,4.5)(0,-0.4)
\put(1,0){\carre}\put(3,0){\carre}
\put(0,1){\carre}\put(2,1){\carre}
\put(1,2){\carre}\put(3,2){\carre}
\put(0,3){\carre}\put(2,3){\carre}
\put(0.5,0.03){\jaune{\circle*{.9}}\hskip-10pt\rouge{\circle{.9}}}
\put(0.25,-0.25){$\rouge{\mathbf{2}}$}
\thicklines
\multiput(0,-0.5)(0,2){3}{\line(1,0){4}}
\multiput(0,-0.5)(2,0){3}{\line(0,1){4}}
\thinlines
\multiput(0,0.5)(0,2){2}{\line(1,0){4}}
\multiput(1,-0.5)(2,0){2}{\line(0,1){4}}
\end{picture}\\ 
           \bf  12 & \bf \bar{1}2 & \bf 1\bar{2} &\bf  \bar{1}\bar{2}\\
\begin{picture}(4.5,4.5)(0,-0.4)
\put(1,0){\carre}\put(3,0){\carre}
\put(0,1){\carre}\put(2,1){\carre}
\put(1,2){\carre}\put(3,2){\carre}
\put(0,3){\carre}\put(2,3){\carre}
\put(1.5,3.03){\jaune{\circle*{.9}}\hskip-10pt\rouge{\circle{.9}}}
\put(1.25,2.75){$\rouge{\mathbf{1}}$}
\put(3.5,1.03){\jaune{\circle*{.9}}\hskip-10pt\rouge{\circle{.9}}}
\put(3.25,.75){$\rouge{\mathbf{1}}$}
\thicklines
\multiput(0,-0.5)(0,2){3}{\line(1,0){4}}
\multiput(0,-0.5)(2,0){3}{\line(0,1){4}}
\thinlines
\multiput(0,0.5)(0,2){2}{\line(1,0){4}}
\multiput(1,-0.5)(2,0){2}{\line(0,1){4}}
\end{picture} &
\begin{picture}(4.5,4.5)(0,-0.4)
\put(1,0){\carre}\put(3,0){\carre}
\put(0,1){\carre}\put(2,1){\carre}
\put(1,2){\carre}\put(3,2){\carre}
\put(0,3){\carre}\put(2,3){\carre}
\put(0.5,2.03){\jaune{\circle*{.9}}\hskip-10pt\rouge{\circle{.9}}}
\put(0.25,1.75){$\rouge{\mathbf{1}}$}
\put(1.5,1.03){\jaune{\circle*{.9}}\hskip-10pt\rouge{\circle{.9}}}
\put(1.25,.75){$\rouge{\mathbf{1}}$}
\thicklines
\multiput(0,-0.5)(0,2){3}{\line(1,0){4}}
\multiput(0,-0.5)(2,0){3}{\line(0,1){4}}
\thinlines
\multiput(0,0.5)(0,2){2}{\line(1,0){4}}
\multiput(1,-0.5)(2,0){2}{\line(0,1){4}}
\end{picture} &
\begin{picture}(4.5,4.5)(0,-0.4)
\put(1,0){\carre}\put(3,0){\carre}
\put(0,1){\carre}\put(2,1){\carre}
\put(1,2){\carre}\put(3,2){\carre}
\put(0,3){\carre}\put(2,3){\carre}
\put(2.5,0.03){\jaune{\circle*{.9}}\hskip-10pt\rouge{\circle{.9}}}
\put(2.25,-.25){$\rouge{\mathbf{1}}$}
\put(1.5,1.03){\jaune{\circle*{.9}}\hskip-10pt\rouge{\circle{.9}}}
\put(1.25,.75){$\rouge{\mathbf{1}}$}
\thicklines
\multiput(0,-0.5)(0,2){3}{\line(1,0){4}}
\multiput(0,-0.5)(2,0){3}{\line(0,1){4}}
\thinlines
\multiput(0,0.5)(0,2){2}{\line(1,0){4}}
\multiput(1,-0.5)(2,0){2}{\line(0,1){4}}
\end{picture} &
\begin{picture}(4.5,4.5)(0,-0.4)
\put(1,0){\carre}\put(3,0){\carre}
\put(0,1){\carre}\put(2,1){\carre}
\put(1,2){\carre}\put(3,2){\carre}
\put(0,3){\carre}\put(2,3){\carre}
\put(2,3){\carre}\put(0.5,2.03){\jaune{\circle*{.9}}\hskip-10pt\rouge{\circle{.9}}}
\put(0.25,1.75){$\rouge{\mathbf{1}}$}
\put(2.5,0.03){\jaune{\circle*{.9}}\hskip-10pt\rouge{\circle{.9}}}
\put(2.25,-.25){$\rouge{\mathbf{1}}$}
\thicklines
\multiput(0,-0.5)(0,2){3}{\line(1,0){4}}
\multiput(0,-0.5)(2,0){3}{\line(0,1){4}}
\thinlines
\multiput(0,0.5)(0,2){2}{\line(1,0){4}}
\multiput(1,-0.5)(2,0){2}{\line(0,1){4}}
\end{picture}\\ 
           \bf  21 &\bf \bar{2}1 &\bf  2\bar{1} &\bf \bar{2}\bar{1}
\end{array}$$\vskip-10pt
\caption{Compact $e$-diagrams for $n=2$.}\label{B2}
\end{figure}
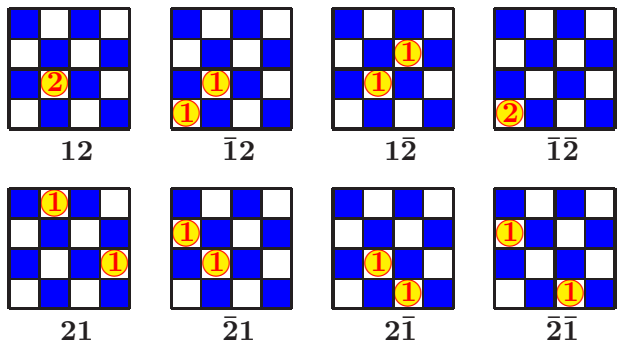

\section{The bijection}

The basic motivation for the introduction of compact $e$-diagram is the following theorem which reflects, in combinatorial term,  the bigraded module isomorphism (\ref{the_iso_inv}), in the case $W=B_n$.

\begin{thm}\label{la_bijection} There is a natural bijection, $\varphi$, between $n$-cell $e$-diagrams and triplets 
   $$D \leftrightarrow (\overline{D},\la,\mu),$$
where $\overline{D}$ is the compactification of $D$, and $\la$ and $\mu$ are two partitions with parts smaller or equal to $n$. Moreover, these partitions are such that 
\begin{equation}\label{prop_bij}
    |D|=|\overline{D}|+2(|\lambda|,|\mu|).
\end{equation}
\end{thm}
\begin{proof} To define $\varphi$, we need only describe the partitions $\lambda=\lambda(D)$ and $\mu=\mu(D)$. This is done as follows. We define the {\em horizontal}, $\mathfrak{h}(D)$, and vertical, $\mathfrak{v}(D)$, {\em marginal distributions} of $D$ as: 
    $$\mathfrak{h}_i(D):=|\{(a,b)\in D\ |\ a>i\}|,$$
and 
   $$\mathfrak{v}_j(D):=|\{(a,b)\in D\ |\ b>i\}|.$$
Considering $\mathfrak{h}(D)$ and $\mathfrak{v}(D)$ as multisets, we then simply define
  $$\lambda(D):=\eta(\mathfrak{h}(D)\setminus \mathfrak{h}(\overline{D})),$$
and
  $$\mu(D):=\eta(\mathfrak{v}(D)\setminus \mathfrak{v}(\overline{D})).$$
Here $\eta$ is the operation that first sorts elements of a multiset in decreasing order, and then erases even indexed entries, as follows:
   $$\eta(a_1,a_2,a_3,a_4,\ldots)=(a_1,a_3,\ldots).$$
As we will currently see, it is a feature of both $\mathfrak{h}(D)\setminus \mathfrak{h}(\overline{D})$ and $\mathfrak{v}(D)\setminus \mathfrak{v}(\overline{D})$ that they are of the form $(a_1,a_2,a_3,a_4,\ldots)$, with 
\begin{equation}\label{parite}
   a_1=a_2,\quad a_3=a_4,\quad  \ldots\ ,
\end{equation}
with $a_i\leq n$.  Both the fact that $\varphi$ is a bijection and property (\ref{parite}) are made evident as follows. The parts of $\lambda(D)$ are simply equal to the number of cells that are moved in each of the following big compacting moves. Using the reading order for cells, we find the smallest cell for which a big left compacting move is possible, and proceed with it. We go on recursively with this process until no left compacting moves remain. Clearly this is reversible if $\lambda$ is known, and a similar description holds for $\mu(D)$. This finishes the proof.
\end{proof}

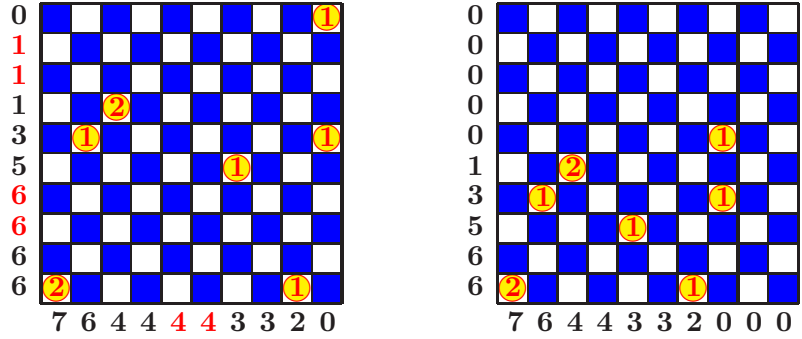
\begin{figure}
\begin{center}
\begin{picture}(11,11)(0,0)
\multiput(1,0)(2,0){5}{\carre}\put(-1,-0.2){${\mathbf 6}$}
\multiput(0,1)(2,0){5}{\carre}\put(-1,.8){${\mathbf 6}$}
\multiput(1,2)(2,0){5}{\carre}\put(-1,1.8){$\rouge{\mathbf 6}$}
\multiput(0,3)(2,0){5}{\carre}\put(-1,2.8){$\rouge{\mathbf 6}$}
\multiput(1,4)(2,0){5}{\carre}\put(-1,3.8){${\mathbf 5}$}
\multiput(0,5)(2,0){5}{\carre}\put(-1,4.8){${\mathbf 3}$}
\multiput(1,6)(2,0){5}{\carre}\put(-1,5.8){${\mathbf 1}$}
\multiput(0,7)(2,0){5}{\carre}\put(-1,6.8){$\rouge{\mathbf 1}$}
\multiput(1,8)(2,0){5}{\carre}\put(-1,7.8){$\rouge{\mathbf 1}$}
\multiput(0,9)(2,0){5}{\carre}\put(-1,8.8){${\mathbf 0}$}
\put(0.5,0.03){\jaune{\circle*{.9}}\hskip-10pt\rouge{\circle{.9}}}\put(0.25,-.25){$\rouge{\mathbf{2}}$}
\put(1.5,5.03){\jaune{\circle*{.9}}\hskip-10pt\rouge{\circle{.9}}}\put(1.25,4.75){$\rouge{\mathbf{1}}$}
\put(2.5,6.03){\jaune{\circle*{.9}}\hskip-10pt\rouge{\circle{.9}}}\put(2.25,5.75){$\rouge{\mathbf{2}}$}
\put(6.5,4.03){\jaune{\circle*{.9}}\hskip-10pt\rouge{\circle{.9}}}\put(6.25,3.75){$\rouge{\mathbf{1}}$}
\put(9.5,5.03){\jaune{\circle*{.9}}\hskip-10pt\rouge{\circle{.9}}}\put(9.25,4.75){$\rouge{\mathbf{1}}$}
\put(8.5,0.03){\jaune{\circle*{.9}}\hskip-10pt\rouge{\circle{.9}}}\put(8.25,-.25){$\rouge{\mathbf{1}}$}
\put(9.5,9.03){\jaune{\circle*{.9}}\hskip-10pt\rouge{\circle{.9}}}\put(9.25,8.75){$\rouge{\mathbf{1}}$}
\thicklines
\multiput(0,-0.5)(0,2){6}{\line(1,0){10}}
\multiput(0,-0.5)(2,0){6}{\line(0,1){10}}
\thinlines
\multiput(0,0.5)(0,2){5}{\line(1,0){10}}
\multiput(1,-0.5)(2,0){5}{\line(0,1){10}}
\put(0,-1.4){$\ \mathbf  7\hskip 5pt \mathbf 6\hskip 5pt\mathbf  4\hskip5pt\mathbf  4\hskip5pt \rouge{\mathbf 4}\hskip5pt \rouge{\mathbf 4}\hskip5pt\mathbf  3\hskip5pt \mathbf 3\hskip5pt \mathbf 2\hskip5pt\mathbf  0$}
\end{picture}
\qquad\qquad
\begin{picture}(11,11)(0,0)
\multiput(1,0)(2,0){5}{\carre}\put(-1,-0.2){${\mathbf 6}$}
\multiput(0,1)(2,0){5}{\carre}\put(-1,.8){${\mathbf 6}$}
\multiput(1,2)(2,0){5}{\carre}\put(-1,1.8){${\mathbf 5}$}
\multiput(0,3)(2,0){5}{\carre}\put(-1,2.8){${\mathbf 3}$}
\multiput(1,4)(2,0){5}{\carre}\put(-1,3.8){${\mathbf 1}$}
\multiput(0,5)(2,0){5}{\carre}\put(-1,4.8){${\mathbf 0}$}
\multiput(1,6)(2,0){5}{\carre}\put(-1,5.8){${\mathbf 0}$}
\multiput(0,7)(2,0){5}{\carre}\put(-1,6.8){${\mathbf 0}$}
\multiput(1,8)(2,0){5}{\carre}\put(-1,7.8){${\mathbf 0}$}
\multiput(0,9)(2,0){5}{\carre}\put(-1,8.8){${\mathbf 0}$}
\put(0.5,0.03){\jaune{\circle*{.9}}\hskip-10pt\rouge{\circle{.9}}}\put(0.25,-.25){$\rouge{\mathbf{2}}$}
\put(1.5,3.03){\jaune{\circle*{.9}}\hskip-10pt\rouge{\circle{.9}}}\put(1.25,2.75){$\rouge{\mathbf{1}}$}
\put(2.5,4.03){\jaune{\circle*{.9}}\hskip-10pt\rouge{\circle{.9}}}\put(2.25,3.75){$\rouge{\mathbf{2}}$}
\put(4.5,2.03){\jaune{\circle*{.9}}\hskip-10pt\rouge{\circle{.9}}}\put(4.25,1.75){$\rouge{\mathbf{1}}$}
\put(7.5,3.03){\jaune{\circle*{.9}}\hskip-10pt\rouge{\circle{.9}}}\put(7.25,2.75){$\rouge{\mathbf{1}}$}
\put(6.5,0.03){\jaune{\circle*{.9}}\hskip-10pt\rouge{\circle{.9}}}\put(6.25,-.25){$\rouge{\mathbf{1}}$}
\put(7.5,5.03){\jaune{\circle*{.9}}\hskip-10pt\rouge{\circle{.9}}}\put(7.25,4.75){$\rouge{\mathbf{1}}$}
\thicklines
\multiput(0,-0.5)(0,2){6}{\line(1,0){10}}
\multiput(0,-0.5)(2,0){6}{\line(0,1){10}}
\thinlines
\multiput(0,0.5)(0,2){5}{\line(1,0){10}}
\multiput(1,-0.5)(2,0){5}{\line(0,1){10}}
\put(0,-1.4){$\ \mathbf  7\hskip 5pt \mathbf 6\hskip 5pt\mathbf  4\hskip5pt\mathbf  4\hskip5pt {\mathbf 3}\hskip5pt {\mathbf 3}\hskip5pt\mathbf  2\hskip5pt \mathbf 0\hskip5pt \mathbf 0\hskip5pt\mathbf  0$}
\end{picture}\vskip10pt
\caption{Differences of marginal partitions of $D$ and $\overline{D}$: $\lambda=4$ and $\mu=61$.}\label{fig_bijection}
\end{center}
\end{figure}

\section{A basis for the trivial component of $\mathcal{C}$}\label{base}

We are ready to describe a straightening algorithm for the expansion of any diagonally invariant polynomial in terms of the elements in the set $\mathcal{M}_n$, with coefficients in the ring $ \QQ[\x]^{B_n} \otimes \QQ[\y]^{B_n}$.
Let $D=({\bf a},{\bf b})$ an $e$-diagram, and consider the effect on $D$ of the bijection $\varphi$ of Theorem \ref{la_bijection}:
  $$\varphi({\bf a},{\bf b})=((\overline{\bf a},\overline{\bf b}),\la,\mu),$$
where we have $(\overline{\bf a},\overline{\bf b})=\overline{D}$.
Then it is not hard to see that
\begin{equation}\label{independence}
M{({\bf a},{\bf b})}=m_{\la}(\x^2)m_{\mu}(\y^2)M(\overline{\bf a},\overline{\bf b}) - \sum_{M' \succ_{lex} M{({\bf a},{\bf b})} } M',
\end{equation}
where $m_{\la}({\bf x}^2)$ is the usual {\em monomial symmetric function} in the squares of the ${\bf x}$ variables. Repeating this process on the remaining terms, we get the desired expansion.

\noindent For example, if $n=2$ and  $ D=\left(\begin{smallmatrix} 1 & 4 \\
                                                     3 & 4  \end{smallmatrix}\right)$; then 
$\overline{D}=\left(\begin{smallmatrix} 1 & 2 \\
                         1 & 2  \end{smallmatrix}\right)$, $\la=1$ and $\mu=11$. 
Then we calculate that 
 $$M{\left( \begin{smallmatrix} 1 & 4 \\
                                   3 & 4  \end{smallmatrix}\right)}=m_1(\x^2)m_{11}(\y^2)M{\left( \begin{smallmatrix} 1 & 2 \\
                                   1 & 2  \end{smallmatrix}\right)}-M{\left( \begin{smallmatrix} 2 & 3 \\
                                   4 & 3  \end{smallmatrix}\right)}.$$
In a similar manner we also get,
 $$M{\left( \begin{smallmatrix} 2 & 3 \\
                                   4 & 3  \end{smallmatrix}\right)}=m_{11}(\x^2)m_{11}(\y^2)M{\left( \begin{smallmatrix} 0 & 1 \\
 2 & 1  \end{smallmatrix}\right)}$$
hence 
$$M{\left( \begin{smallmatrix} 1 & 4 \\
                                   3 & 4  \end{smallmatrix}\right)}=m_1(\x^2)m_{11}(\y^2)M_{1\overline{2}}-m_{11}(\x^2)m_{11}(\y^2)M_{2\overline{1}}.$$

It follows that the $M_\beta$'s are indeed a set of generators for the trivial component of $\mathcal{C}$. Since the dimension of $\mathcal{C}^{B_n}$ is $|B_n|$, as we pointed out at the beginning of Section \ref{triv-section}, we have  
\begin{prop} The set 
  $$\{M_\beta+ \mathcal{I}_{B_n \times B_n} \ |\ \be \in B_n \},$$
is a bihomogeneous basis for the trivial component of $\mathcal{C}$. \qed
\end{prop} 
In particular, the recursive procedure (\ref{independence}) give us also expressions for the invariant polynomials $u_\beta(\x,\y)$ in (\ref{decomposition_unique}). It follows that (\ref{hilbert_compact}) is indeed the bigraded Hilbert series of both $\mathcal{C}^{B_n}$ and $\mathcal{H}^{B_n}$.

\section{The alternating component of $\mathcal{C}$}

Recall that a polynomial
$p(\x,\y)$ is said to $B_n$-diagonally alternating 
  $$ \be \cdot p(\x,\y)={\rm sign} (\be) p(\x,\y),$$
for all $\be \in B_n$. Recall also that the {\em sign} of an element $\beta$ of $B_n$ can be computed as
   $${\rm sign}(\beta)=\frac{\beta\cdot \Delta(\x)}{\Delta(\x)},$$
with $\Delta(\x)$ as given in (\ref{jacobian}).
It is easy to see that a basis of the space $R^{\pm}$ of diagonally $B_n$-alternating polynomials is obtained as follows. An $o$-diagram (``$o$'' for odd)
\begin{equation}\label{o_diag}
    D= ({\bf a},{\bf b})= \begin{pmatrix} a_1 & a_2 & \ldots & a_n\\ 
                       b_1 & b_2 & \ldots & b_n \end{pmatrix}
\end{equation}
is any $n$-element subset of $\NN \times \NN$, with all cells of {\em odd parity}. This is to say all the $a_i+b_i$'s are odd. As is now our custom, the cells of $D$ are ordered lexicographically in (\ref{o_diag}). We emphasize that in the present context the cells of an $o$-diagram are all distinct. We then consider the determinant
\begin{equation}
\Delta_D(\x,\y):=\begin{array}{|cccc|}
x_1^{a_1}y_1^{b_1} & x_1^{a_2}y_1^{b_2} & \ldots & x_1^{a_n}y_1^{b_n} \\
\vdots & \vdots & \vdots & \vdots \\
x_n^{a_1}y_n^{b_1} & x_n^{a_2}y_n^{b_2} & \ldots & x_n^{a_n}y_n^{b_n} 
\end{array}\,,
\end{equation}
which is easily seen to be  diagonally $B_n$-alternating.  It is not hard to see that a linear basis for $R^{\pm}$ is given by the set 
\begin{equation}\label{base_alt}
    \{\Delta_D  \ |\ D \subseteq (\NN \times \NN)^1, \ |D|=n\  \}.
\end{equation}
Here $(\NN \times \NN)^1$ is the  {\em odd-chessboard plane}:
 $$ (\NN \times \NN)^1:=\{(a,b) \in \NN \times \NN: a+b \equiv 1 \ ({\rm mod}\ 2)\}.$$
 Thus an $n$ element subset of $(\NN \times \NN)^1$ is just an $o$-diagram.
Using (\ref{the_iso}), (\ref{frob_R}), and the fact that the {\em sign} representation of $B_n$ has Frobenius equal to $e_n[\z-\bar{\z}]$, we get that
   $$e_n\left[\frac{q+t}{(1-q^2)(1-t^2)}\right]$$
is the bigraded Hilbert series of $R^\pm$.

As we have seen previously, the bijection $\varphi$ of Theorem \ref{la_bijection} is a combinatorial ``shadow'' of the bigraded module isomorphism (\ref{the_iso_inv}). In the same manner, we can translate in combinatorial terms the isomorphism (\ref{the_iso_alt}). This involves a similar bijection, but a different notion of compact diagrams. In fact, these new compact diagrams naturally appear as transformations of compact $e$-diagrams. This goes through a combinatorial  ``interpretation'' of the following linear operators. 

\begin{lem}\label{inv_to_alt}
  As linear operators on $\mathcal{H}$, both
        $$p({\bf x},{\bf y})\mapsto p(\partial {\bf x},{\bf y})\Delta({\bf x}),\qquad {\rm and}\qquad
            p({\bf x},{\bf y})\mapsto p( {\bf x},\partial {\bf y})\Delta({\bf y}),$$
send diagonally $B_n$-invariant polynomials to diagonally  $B_n$-alternating polynomials.
\end{lem}
We can mimic the effect of these operators with a map that sends compact $e$-diagrams to some new special $o$-diagrams:
\begin{equation}\label{flipx}
  \psi:({\bf g}(\beta),\widetilde{{\bf g}}(\beta))\mapsto ({\bf g}(\beta),{\bf c}(\beta)-\widetilde{{\bf g}}(\beta)),
\end{equation}
where ${\bf c}(\beta)=(c_1,\ldots,c_n)$ is the permutation of the integers $1,3,\ldots,(2n-1)$ obtained as
   $${\bf c}(\beta)=(2\sigma_1-1,2\sigma_2-1, \ldots, 2\sigma_n-1),$$
 with $\sigma_i:=|\beta(i)|$. 
It is clear that the cells of $\psi(D_\beta)$  are all of odd parity. As will be checked below, the image under $\psi$ of a compact $e$-diagram is always an $o$-diagram. For instance, the image:
\begin{equation}\label{un_o_diag}
    \begin{pmatrix} 0 & 0 & 1  & 2 & 2 & 4 & 6 & 7 & 7 \\
                             1 & 3 & 6  & 9 & 11 & 5 & 5 & 8 & 12 \end{pmatrix},
\end{equation}                             
under $\psi$,  of the compact $e$-diagram of (\ref{un_diag_comp}), is illustrated in Figure \ref{fig_compact}.
\setlength{\unitlength}{3mm}
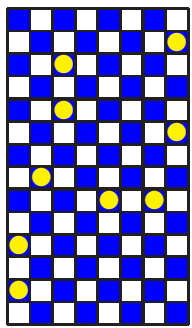
\begin{figure}\begin{center}
 \begin{picture}(11,13)(0,0)
\multiput(1,0)(2,0){4}{\carre}
\multiput(0,1)(2,0){4}{\carre}
\multiput(1,2)(2,0){4}{\carre}
\multiput(0,3)(2,0){4}{\carre}
\multiput(1,4)(2,0){4}{\carre}
\multiput(0,5)(2,0){4}{\carre}
\multiput(1,6)(2,0){4}{\carre}
\multiput(0,7)(2,0){4}{\carre}
\multiput(1,8)(2,0){4}{\carre}
\multiput(0,9)(2,0){4}{\carre}
\multiput(1,10)(2,0){4}{\carre}
\multiput(0,11)(2,0){4}{\carre}
\multiput(1,12)(2,0){4}{\carre}
\multiput(0,13)(2,0){4}{\carre}
\thicklines
\multiput(0,-0.5)(0,2){8}{\line(1,0){8}}
\multiput(0,-0.5)(2,0){5}{\line(0,1){14}}
\thinlines
\multiput(0,0.5)(0,2){7}{\line(1,0){8}}
\multiput(1,-0.5)(2,0){4}{\line(0,1){14}}
\put(0.5,1.03){\jaune{\circle*{.8}}}
\put(0.5,3.03){\jaune{\circle*{.8}}}
\put(1.5,6.03){\jaune{\circle*{.8}}}
\put(2.5,9.03){\jaune{\circle*{.8}}}
\put(2.5,11.03){\jaune{\circle*{.8}}}
\put(4.5,5.03){\jaune{\circle*{.8}}}
\put(6.5,5.03){\jaune{\circle*{.8}}}
\put(7.5,8.03){\jaune{\circle*{.8}}}
\put(7.5,12.03){\jaune{\circle*{.8}}}
\end{picture}
\caption{Image under $\psi$ of a compact $e$-diagram.}\label{fig_compact}
\end{center}\end{figure}
 
 Just as in the $e$-diagram case, there is a classification of $o$-diagrams in terms of elements of $B_n$. The construction is very similar, with only a change in the labeling order. In this context, we rather use the increasing {\em colabeling} order, ``$\prec_{col}$'', for which
 $$(a,b)\prec_{col} (a',b') \iff \begin{cases}
      b < b'  & \text{or}, \\
      b=b' & \text{and}\quad a>a'.
\end{cases}$$
Thus we get
\begin{equation}
    \begin{pmatrix} 0 & 0 & 6 & 4 & 1 & 7 & 2 & 2  & 7 \\
                             1 & 3 & 5 & 5 & 6 & 8  & 9 & 11 & 12 \end{pmatrix},
\end{equation}                             
as the colabeling order sorting of the diagram in (\ref{un_o_diag}). The {\em classifying signed permutation} of an $o$-diagram is then obtained just as  in Section \ref{combin}, starting with this new labeling order. Just as in our previous case, we say that two $o$-diagrams are equivalent $D\simeq \widetilde{D}$, if they have the same classifying signed permutation. The reader may check that $\bar{1}\bar{2}5 \bar{7}\bar{8} \bar{4}\bar{3} 6 9$ is the classifying signed permutation of the $o$-diagram in Figure \ref{fig_compact}.

\section{Compact $o$-diagrams}
Going on with an approach analogous to that of Section \ref{sec_comp}, we construct a minimal $o$-diagram among those classified by a given $\beta$ in $B_n$.  The {\em compact $o$-diagram} associated to $\beta$ is
\begin{equation}\label{defodiagram}
 D_{\be}^s:=\begin{pmatrix} g_1(\be) & g_2(\be) & \ldots & g_n(\be)\\
                            \widehat{g}_{1}(\be) & \widehat{g}_{2}(\be) & \ldots & \widehat{g}_{n}(\be) \end{pmatrix},
\end{equation}
where $g_i(\beta)$ is defined as in (\ref{def-g_i}) and \[\widehat{g}_i(\be):=2\mu_{\sigma(i)}(\beta)+\varepsilon_i(\be).\]
Here $\varepsilon_i(\beta)$ is defined as in (\ref{varepsilon}), $\sigma(i)$ is the absolute value of $\beta(i)$, and
   $$\mu_i(\be):=|\{k \in \Ris(\be^{-1}): k <i\},$$
where the set $\Ris(\be)$ is the set of {\em rises} of $\be$, namely
    $$\Ris(\be):=\{i \in [n-1]: \be(i)<\be(i+1)\}.$$
For $n=2$, the   8 compact  $o$-diagrams are illustrated in Figure \ref{liste_o_comp}.
\setlength{\unitlength}{4mm}
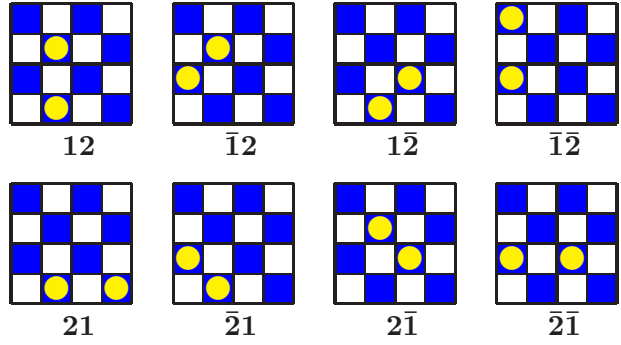
\begin{figure}
 $$\begin{array}{cccc}
\begin{picture}(4.5,4.5)(0,-0.4)
\put(1,0){\carre}\put(3,0){\carre}
\put(0,1){\carre}\put(2,1){\carre}
\put(1,2){\carre}\put(3,2){\carre}
\put(0,3){\carre}\put(2,3){\carre}
\put(1.5,2.03){\jaune{\circle*{.8}}}
\put(1.5,0.03){\jaune{\circle*{.8}}}
\thicklines
\multiput(0,-0.5)(0,2){3}{\line(1,0){4}}
\multiput(0,-0.5)(2,0){3}{\line(0,1){4}}
\thinlines
\multiput(0,0.5)(0,2){2}{\line(1,0){4}}
\multiput(1,-0.5)(2,0){2}{\line(0,1){4}}
\end{picture}&
\begin{picture}(4.5,4.5)(0,-0.4)
\put(1,0){\carre}\put(3,0){\carre}
\put(0,1){\carre}\put(2,1){\carre}
\put(1,2){\carre}\put(3,2){\carre}
\put(0,3){\carre}\put(2,3){\carre}
\put(1.5,2.03){\jaune{\circle*{.8}}}
\put(0.5,1.03){\jaune{\circle*{.8}}}
\thicklines
\multiput(0,-0.5)(0,2){3}{\line(1,0){4}}
\multiput(0,-0.5)(2,0){3}{\line(0,1){4}}
\thinlines
\multiput(0,0.5)(0,2){2}{\line(1,0){4}}
\multiput(1,-0.5)(2,0){2}{\line(0,1){4}}
\end{picture} &
\begin{picture}(4.5,4.5)(0,-0.4)
\put(1,0){\carre}\put(3,0){\carre}
\put(0,1){\carre}\put(2,1){\carre}
\put(1,2){\carre}\put(3,2){\carre}
\put(0,3){\carre}\put(2,3){\carre}
\put(2.5,1.03){\jaune{\circle*{.8}}}
\put(1.5,0.03){\jaune{\circle*{.8}}}
\thicklines
\multiput(0,-0.5)(0,2){3}{\line(1,0){4}}
\multiput(0,-0.5)(2,0){3}{\line(0,1){4}}
\thinlines
\multiput(0,0.5)(0,2){2}{\line(1,0){4}}
\multiput(1,-0.5)(2,0){2}{\line(0,1){4}}
\end{picture} &
\begin{picture}(4.5,4.5)(0,-0.4)
\put(1,0){\carre}\put(3,0){\carre}
\put(0,1){\carre}\put(2,1){\carre}
\put(1,2){\carre}\put(3,2){\carre}
\put(0,3){\carre}\put(2,3){\carre}
\put(0.5,3.03){\jaune{\circle*{.8}}}
\put(0.5,1.03){\jaune{\circle*{.8}}}
\thicklines
\multiput(0,-0.5)(0,2){3}{\line(1,0){4}}
\multiput(0,-0.5)(2,0){3}{\line(0,1){4}}
\thinlines
\multiput(0,0.5)(0,2){2}{\line(1,0){4}}
\multiput(1,-0.5)(2,0){2}{\line(0,1){4}}
\end{picture}\\ 
           \bf  12 & \bf \bar{1}2 & \bf 1\bar{2} &\bf  \bar{1}\bar{2}\\
 \begin{picture}(4.5,4.5)(0,-0.4)
\put(1,0){\carre}\put(3,0){\carre}
\put(0,1){\carre}\put(2,1){\carre}
\put(1,2){\carre}\put(3,2){\carre}
\put(0,3){\carre}\put(2,3){\carre}
\put(1.5,0.03){\jaune{\circle*{.8}}}
\put(3.5,0.03){\jaune{\circle*{.8}}}
\thicklines
\multiput(0,-0.5)(0,2){3}{\line(1,0){4}}
\multiput(0,-0.5)(2,0){3}{\line(0,1){4}}
\thinlines
\multiput(0,0.5)(0,2){2}{\line(1,0){4}}
\multiput(1,-0.5)(2,0){2}{\line(0,1){4}}
\end{picture}&
\begin{picture}(4.5,4.5)(0,-0.4)
\put(1,0){\carre}\put(3,0){\carre}
\put(0,1){\carre}\put(2,1){\carre}
\put(1,2){\carre}\put(3,2){\carre}
\put(0,3){\carre}\put(2,3){\carre}
\put(1.5,0.03){\jaune{\circle*{.8}}}
\put(0.5,1.03){\jaune{\circle*{.8}}}
\thicklines
\multiput(0,-0.5)(0,2){3}{\line(1,0){4}}
\multiput(0,-0.5)(2,0){3}{\line(0,1){4}}
\thinlines
\multiput(0,0.5)(0,2){2}{\line(1,0){4}}
\multiput(1,-0.5)(2,0){2}{\line(0,1){4}}
\end{picture}&
\begin{picture}(4.5,4.5)(0,-0.4)
\put(1,0){\carre}\put(3,0){\carre}
\put(0,1){\carre}\put(2,1){\carre}
\put(1,2){\carre}\put(3,2){\carre}
\put(0,3){\carre}\put(2,3){\carre}
\put(1.5,2.03){\jaune{\circle*{.8}}}
\put(2.5,1.03){\jaune{\circle*{.8}}}
\thicklines
\multiput(0,-0.5)(0,2){3}{\line(1,0){4}}
\multiput(0,-0.5)(2,0){3}{\line(0,1){4}}
\thinlines
\multiput(0,0.5)(0,2){2}{\line(1,0){4}}
\multiput(1,-0.5)(2,0){2}{\line(0,1){4}}
\end{picture}  &
\begin{picture}(4.5,4.5)(0,-.4)
\put(1,0){\carre}\put(3,0){\carre}
\put(0,1){\carre}\put(2,1){\carre}
\put(1,2){\carre}\put(3,2){\carre}
\put(0,3){\carre}\put(2,3){\carre}
\put(0.5,1.03){\jaune{\circle*{.8}}}
\put(2.5,1.03){\jaune{\circle*{.8}}}
\thicklines
\multiput(0,-0.5)(0,2){3}{\line(1,0){4}}
\multiput(0,-0.5)(2,0){3}{\line(0,1){4}}
\thinlines
\multiput(0,0.5)(0,2){2}{\line(1,0){4}}
\multiput(1,-0.5)(2,0){2}{\line(0,1){4}}
\end{picture} \\ 
           \bf  21 &\bf \bar{2}1 &\bf  2\bar{1} &\bf \bar{2}\bar{1}
\end{array}$$
\caption{All 8 compact $o$-diagrams with 2 cells.}\label{liste_o_comp}
\end{figure}
The point of all this is that
\begin{equation}\label{point}
    D_\beta^s=\psi(D_\beta).
\end{equation}
We can thus translate results on compact $e$-diagram into results on $o$-compact diagrams. For instance, we get
\begin{cor}
Let $n \in \NN$. Then
   $$\sum_{\be \in B_n}q^{|{\bf g}(\be)|}t^{|\widehat{\bf g}(\be)|}=\sum_{\be \in B_n}q^{\fmaj(\be)}t^{n^2-\fmaj(\be^{-1})}.$$
\end{cor}
Along the same lines, we defined the {\em compactification} of $o$-diagrams as
 \begin{equation}\label{o_comp}
        \overline{D}^s:=D_\beta^s,
  \end{equation}
with $\beta$ being the classifying signed permutation for the $o$-diagram $D$.
Following an argument close to that of Theorem \ref{la_bijection}, we get a bijection
   $$D\leftrightarrow (\overline{D}^s,\lambda,\mu)$$
with properties as before. Hence,
\begin{equation}
H_{q,t}(R^\pm)= e_n\left[\frac{q+t}{(1-q^2)(1-t^2)}\right]=\sum_{\be \in B_n}q^{\fmaj(\be)}t^{n^2-\fmaj(\be^{-1})}
 \prod_{i=1}^n \frac{1}{1-q^{2i}}\prod_{i=1}^n \frac{1}{1-t^{2i}}.\label{ogenfunction}
\end{equation}
There are similar identities associated to each irreducible character of $B_n$, taking the form
\begin{equation*}
s_\lambda\left[\frac{1+q\,t}{(1-q^2)(1-t^2)}\right] s_\mu\left[\frac{q+t}{(1-q^2)(1-t^2)}\right]=\left[n\atop |\lambda|\right]_{q^2}\left[n\atop |\mu|\right]_{t^2} \Psi_{\lambda,\mu}(q,t)
 \prod_{i=1}^n \frac{1}{1-q^{2i}}\prod_{i=1}^n \frac{1}{1-t^{2i}}
\end{equation*}  
with expressions in bracket (in the right hand side) standing for $q^2$-binomial coefficients (or $t^2$-binomial coefficients), and $\Psi_{\lambda,\mu}(q,t)$ a positive integer coefficient polynomial. 
These identities can also be explained through bijections similar to those that we have considered above.

\end{document}